\documentclass[11pt]{article}
\title{Pairing and duality of algebraic quantum groups}
\author{Thomas Timmermann, Alfons Van Daele and Shuanhong Wang}
\usepackage{amsmath}
\usepackage{amsthm}
\usepackage{amssymb}
\usepackage{hyperref}
\usepackage{color}
\newtheorem{thm}{Theorem}
\newtheorem{inspr}[thm]{}
\newenvironment{definitie}{\begin{itemize}\item[ ]\hspace{-26pt}\bf Definition \rm }{\end{itemize}}
\newenvironment{notatie}{\begin{itemize}\item[ ]\hspace{-26pt}\bf Notation \rm }{\end{itemize}}
\newenvironment{voorbeeld}{\begin{itemize}\item[ ]\hspace{-26pt}\bf Example \rm }{\end{itemize}}
\newenvironment{stelling}{\begin{itemize}\item[ ]\hspace{-26pt}\bf Theorem \rm }{\end{itemize}}
\newenvironment{propositie}{\begin{itemize}\item[ ]\hspace{-26pt}\bf Proposition \rm }{\end{itemize}}
\newenvironment{lemma}{\begin{itemize}\item[ ]\hspace{-26pt}\bf Lemma \rm }{\end{itemize}}
\newenvironment{opmerking}{\begin{itemize}\item[ ]\hspace{-26pt}\bf Remark \rm }{\end{itemize}}
\newenvironment{voorwaarde}{\begin{itemize}\item[ ]\hspace{-26pt}\bf Assumption \rm }{\end{itemize}}

\newcommand{\defin}{\begin{inspr}\begin{definitie}}
\newcommand{\edefin}{\end{definitie}\end{inspr}}

\newcommand{\notat}{\begin{inspr}\begin{notatie}}
\newcommand{\enotat}{\end{notatie}\end{inspr}}

\newcommand{\voorb}{\begin{inspr}\begin{voorbeeld}}
\newcommand{\evoorb}{\end{voorbeeld}\end{inspr}}

\newcommand{\stel}{\begin{inspr}\begin{stelling}}
\newcommand{\estel}{\end{stelling}\end{inspr}}

\newcommand{\prop}{\begin{inspr}\begin{propositie}}
\newcommand{\eprop}{\end{propositie}\end{inspr}}

\newcommand{\lem}{\begin{inspr}\begin{lemma}}
\newcommand{\elem}{\end{lemma}\end{inspr}}

\newcommand{\opm}{\begin{inspr}\begin{opmerking}}
\newcommand{\eopm}{\end{opmerking}\end{inspr}}

\newcommand{\voorw}{\begin{inspr}\begin{voorwaarde}}
\newcommand{\evoorw}{\end{voorwaarde}\end{inspr}}

\newcommand{\bew}{\vspace{-0.3cm}\begin{itemize}\item[ ] \bf Proof\rm: }
\newcommand{\ebew}{\hfill $\qed$ \end{itemize}}

\newcommand{\ssnl}{\vskip 3pt}

\newcommand{\snl}{\vskip 7pt}
\newcommand{\nl}{\vskip 12pt}

\newcommand{\inv}{^{-1}}
\newcommand{\ot}{\otimes}

\newcommand{\tl}{\triangleleft}
\newcommand{\tr}{\triangleright}

\newcommand{\tussenen}{\qquad\quad\text{and}\qquad\quad}

\numberwithin{thm}{section}
\numberwithin{equation}{section}

\begin{document}

\centerline{\bf \Large Pairing and duality of algebraic quantum groupoids}
\vspace{13 pt}
\centerline{\it T.\ Timmermann  $^{(1)}$, A.\ Van Daele \rm $^{(2)}$ \rm and \it S.H.\ Wang \rm $^{(3)}$}
\vspace{20 pt}
{\bf Abstract}
\nl
\emph{Algebraic quantum groupoids} have been developed by two of the authors of this note (AVD and SHW) in a series of papers \cite{VD-W3, VD-W4, VD-W5} and \cite{VD-W6}, see also \cite{VD8}. By an algebraic quantum groupoid, we understand a \emph{regular weak multiplier Hopf algebra with enough integrals}. Regular multiplier \emph{Hopf algebroids} are obtained also by two authors of this note (TT and AVD) in \cite{T-VD1}. Integral theory and duality for those have been studied by one author here (TT) in \cite{T1,T2}. In these papers, the term \emph{algebraic quantum groupoid} is used for a regular multiplier Hopf \emph{algebroid} with a \emph{single} faithful integral. Finally, again two authors of us (TT and AVD) have investigated the relation between weak multiplier Hopf algebras and multiplier Hopf algebroids in \cite{T-VD2}.
\ssnl
In the paper \emph{Weak multiplier Hopf algebras III. Integrals and duality} \cite{VD-W6}, one of the main results is that the dual of an algebraic quantum groupoid, admits a dual of the same type. In the paper \emph{On duality of algebraic quantum groupoids} \cite{T2}, a result of the same nature is obtained for regular multiplier Hopf algebroids with a single faithful integral. The duality of regular weak multiplier Hopf algebras with a single integral can be obtained from the duality of regular multiplier Hopf algebroids (see \cite{T2}). That is however  not the obvious way to obtain this result. It is more difficult and less natural than the direct way followed in \cite{VD-W6}. We will discuss this statement further in the paper.
\ssnl
Nevertheless, it is interesting to investigate the relation between the two approaches to duality in greater detail. This is what we do in this paper. We build further on the intimate relation between weak multiplier Hopf algebras and multiplier Hopf algebroids as studied in \cite{T-VD2}. We now add the presence of integrals. That seems to be done best in a framework of dual pairs. It is in fact more general than the duality of these objects coming with integrals.
\ssnl
We are convinced that the material we present in this paper will provide a deeper understanding of the duality of algebraic quantum groupoids, both within the framework of weak multiplier Hopf algebras, as well as more generally for multiplier Hopf algebroids.

\ssnl
Finally, we feel it is also appropriate to include some historical comments on the development of these duality theories.
\nl
\vskip 15 pt
\hrule
\vskip 7 pt
\begin{itemize}
\item[(1)] {\it Address}: Enschedeweg 63, D-48149 M\"unster, Germany.\newline {\it E-mail}: \texttt{thomas@timmer-net.de}
\item[(2)] {\it Address}: Department of Mathematics, University of Leuven, Celestijnenlaan 200B,\newline
B-3001 Heverlee (Belgium). {\it Email}: \texttt{alfons.vanDaele@kuleuven.be}
\item[(3)] {\it  Corresponding author's address}: Shing-Tung Yau Center, School of Mathematics, Southeast University, Nanjing 210096, China.
{\it Email}: \texttt{shuanhwang@seu.edu.cn}
\end{itemize}

\newpage

\setcounter{section}{-1}  

\section{\hspace{-17pt}. Introduction}   \label{s:introduction} 

A  weak multiplier Hopf algebra is a pair $(A,\Delta)$ of an algebra $A$ and a coproduct $\Delta$ on $A$. The algebra is not necessarily unital but the product is assumed to be non-degenerate. If $A$ has a unit, the coproduct is not required to be unital.
More generally, if there is no unit in $A$, the coproduct is required to be non-degenerate. Recall that for a non-degenerate coproduct, there exists a unital extension to the multiplier algebra. There are more conditions on $\Delta$. They are  such that the function algebra on a groupoid $G$, with the coproduct induced by the multiplication, satisfies the requirements. The function algebra is the space of complex functions on $G$ with finite support with pointwise product and the coproduct is induced by the product in $G$.
\ssnl
For details we refer to our work on weak multiplier Hopf algebras, see under the item \emph{Basic references}. For the convenience of the reader, we have included a summary of the most important definitions and necessary results in Section \ref{s:prel} with more preliminaries. The source and target algebra play an important role. For the function algebra on a groupoid, they are constructed with the source and target maps. Again in Section \ref{s:prel}, we will spend a good part to recall the definitions and properties.
\ssnl
To any regular weak multiplier Hopf algebra $A$ is canonically associated a regular multiplier Hopf algebroid. The process is essentially based on the passage from the tensor $A\ot A$ of $A$ with itself to appropriate balanced tensor products. These are constructed using actions coming from multiplication, left or right, with the source and target algebras. Details are found in \cite{T-VD2} but also here we will recall the basic aspects in the preliminary Section \ref{s:prel}. We recall the notion of a regular multiplier Hopf algebroid while discussing the passage from weak multiplier Hopf algebras to algebroids.
\ssnl
As mentioned already in the abstract, in this note we use the term \emph{algebraic quantum groupoid} for a regular weak multiplier Hopf algebra $(A,\Delta)$ with enough integrals.
The notion will be recalled in Section \ref{s:from}. There we will see how integrals on weak multiplier Hopf algebras give rise to  integrals on the associated multiplier Hopf algebroid. This is the main content of Section \ref{s:from}. The method is slightly different from how it is done in Section 5.1 of \cite{T2}. In a way it is a more direct approach from what is found there. In particular, we avoid the use of results from \cite{T1}.
\ssnl
In Section \ref{s:dual1} we start the study of dual pairs of weak multiplier Hopf algebras $A$ and $A'$. The pair $(A,\widehat A)$, where $A$ is a regular multiplier Hopf algebra with enough integrals and $\widehat A$ its dual, is an example of such a pairing. In the general setting, we see how the pairing behaves with respect to the source and target algebras of both $A$ and $A'$. This behavior is crucial in order to associate the pairings on balanced tensor products.
We will do this in Section \ref{s:dual2} where we show how this case gives rise to a pair of multiplier Hopf algebroids.
\ssnl
Indeed, in Section \ref{s:dual2} we consider dual pairs of multiplier Hopf algebroids $\mathfrak A$ and $\mathfrak A'$. We begin with imposing some natural conditions so that we can give a definition of such a pair. It should be considered rather as a preliminary approach to the study of dual pairs of multiplier Hopf algebroids.
\ssnl
We consider two special cases that comply with this preliminary notion in Section \ref{s:dual3}. First we have the dual pair that arises from a dual pair of weak multiplier Hopf algebras by applying the procedure of Section \ref{s:prel}. We associate  with $A$ and $A'$ the algebroids as in Section \ref{s:prel} and we show that this gives rise to a pair of multiplier Hopf algebroids. A second case is where $\mathfrak A$ is a multiplier Hopf algebroid with a single faithful integral and $\mathfrak A'$ is the dual in the sense duality of multiplier Hopf algebras with integrals as in \cite{T2} . When we start from a regular weak multiplier Hopf algebra with a single faithful integral, we can perform these two steps in a different order. On the one hand, we can define the dual in the sense of duality of algebraic quantum groupoids as in \cite{VD-W6} and then pass to the pairing of the associated multiplier Hopf algebroids. On the other hand, we can first associate the multiplier Hopf algebroid and then take the dual in the sense of \cite{T2}. We show that the resulting pair of algebroids is the same. These examples are important for a better understanding of the whole theory.
\ssnl
In Section 5.1 of \cite{T2}), in the case of a \emph{single faithful integral}, the dual of an algebraic quantum groupoid is obtained as a special case of the duality of regular multiplier Hopf algebroids with a single faithful integral (see Section 5.1 of \cite{T2}). This is of course intimately related with the result we prove in our section on duality in this paper as mentioned above.
\ssnl
Indeed, as we understand from the previous observations, there are two ways to obtain the dual of a weak multiplier Hopf algebra with integrals. There is the direct way, along the lines of previous work on weak multiplier Hopf algebras \cite{VD-W3, VD-W4, VD-W5}, as in \cite{VD-W6} and there is the possibility to see it as a consequence of the duality for regular multiplier Hopf algebroids with integrals as it is done in \cite{T2}. We will comment in great detail on this relation and how they compare with each other.
\ssnl
Moreover, as already mentioned in the abstract, we will add to this discussion some historical considerations about how these results were developed, in relation to each other. This is done in Section \ref{s:hist}. We also consider this aspect as an important contribution to the discussion in this note.
\ssnl
In the last section, Section \ref{s:concl}, we finish with conclusions and add some more ideas for future research on this material.
\nl
\bf The place of these results in the theory \rm
\nl
What is our intention with writing this note?
\ssnl
First remark that the concept of a multiplier Hopf algebroid is not an easy one. This is even more so for the theory of integrals on multiplier Hopf algebroids. Moreover the notion of a measured multiplier Hopf algebroid (a regular multiplier Hopf algebroid with a faithful integral), as it is defined in Definition 2.2.1 in \cite{T2}, involves some unexpected conditions (see e.g.\ condition (4) in that definition). Finally, remark that the duality of regular weak multiplier Hopf algebras as obtained in Section 5.1 of \cite{T2} is more involved than it looks at a first glance.
\ssnl
One of the points we make in this note is that the direct approach to the duality of regular weak multiplier Hopf algebras is after all easier, more natural, and in fact also more general. It is an obvious step in the development of the theory of weak multiplier Hopf algebras. We also have more results on integrals to begin with in \cite{VD-W6}. Remark further that integrals on weak multiplier Hopf algebras appear already in the paper on the Larson Sweedler theorem \cite{K-VD}.
\ssnl
Another point we make with this note is that first studying the easier theory of duality for weak multiplier Hopf algebras will help to understand the more general, but also more complicated duality theory for multiplier Hopf algebroids. The framework of pairings that we use here is more general and still, in some sense, easier.
\ssnl
We include a fair amount of details in this note. Precisely because we also want to contribute to a better understanding of duality for multiplier Hopf algebroids with integrals.
\nl
\bf Notations and conventions \rm
\nl
We consider  algebras $A$ over the field of complex numbers. They need not have a unit but the product is always assumed to be non-degenerate. Then we have the multiplier algebra $M(A)$ of $A$. It is the largest algebra with identity that contains $A$ as an essential two-sided ideal. We use $A^\text{op}$ for the algebra obtained from $A$ by reversing the product.
\ssnl
A coproduct $\Delta$ on an algebra $A$ is a homomorphism from $A$ to the multiplier algebra $M(A\ot A)$ satisfying coassociativity. However, in order to express coassociativity for such a map, some regularity conditions are needed. This is documented in the literature (see the basic works on multiplier Hopf algebras and weak multiplier Hopf algebras). If $\Delta$ is a coproduct on $A$ we use $\Delta^\text{cop}$ for the coproduct on $A$ obtained by composing $\Delta$ with the flip map.
\ssnl
The coproduct is not assumed to be unital (in the case of a unital algebra) and more generally, it is not assumed to be non-degenerate. Some weaker condition is assumed. This still guarantees that the coproduct can be extended to the multiplier algebra, but also this extension is not necessarily unital. See the original papers on weak multiplier Hopf algebras for these properties.
\ssnl
We also sometimes use the Sweedler notation for the coproduct. This has to be done with care by arguing that the necessary coverings can be provided. The use of the Sweedler notation for coproducts with values in the multiplier algebra has been documented in various papers on the subject. See e.g. \cite{Dr-VD} and \cite{VD2a}. A new note with more information about the use of the Sweedler notation in this context is being prepared, see \cite{VD7}.
\ssnl
We use $\iota$ to denote the identity map on various spaces.
\ssnl
Remark finally that notations used in the papers on multiplier Hopf algebroids differ somewhat from those used in the theory of weak multiplier Hopf algebras. For the most part, we will use the notations as in the papers on multiplier Hopf algebroids \cite{T-VD1, T1, T2}. These notations are better as they cover more information.
On the other hand we do avoid the use of too many subscripts and superscripts.
\ssnl
Further,  we partly have to rely on the material of \cite{T-VD2} where the notations from weak multiplier Hopf algebras are used. We are aware of the confusion this may cause and for this reason, where appropriate, we will provide a \emph{dictionary} for the notations used in the two sets of papers where appropriate.
\nl
\bf Basic references \rm
\nl
For the theory of Hopf algebras, we refer to the original work by Abe \cite{A} and Sweedler \cite{S}, as well as to the more recent treatment by Radford \cite{Ra}. Multiplier Hopf algebras are studied in \cite{VD1} and multiplier Hopf algebras with integrals (algebraic quantum groups) in \cite{VD2}. The theory of weak Hopf algebras is developed in \cite{B-N-S} and \cite{B-S} but also the works on finite quantum groupoids in \cite{N-V1} and \cite{N-V2} are very useful. Weak multiplier Hopf algebras are developed in a series of papers \cite{VD-W3, VD-W4, VD-W5} and \cite{VD-W6}. In particular \cite{VD-W5} depends on the theory of separability idempotents as studied in \cite{VD4.v1} (in the regular case) and \cite{VD4.v2} (in the more general case). For multiplier Hopf algebroids, the references are \cite{T-VD1} for the basic theory, \cite{T-VD2} for the relation with weak multiplier Hopf algebras and \cite{T1,T2} for integrals and duality of multiplier Hopf algebroids.
\nl
\bf Acknowledgments  \rm
\nl
 The authors are very grateful to the anonymous referee for his/her thorough review of this work and his/her comments.
\ssnl
The first-named author (Thomas Timmermann) would like to thank the second-named author (Alfons Van Daele) for his on-going support over the last years. The last name author (Shuanhong Wang) would like to thank his coauthor Alfons Van Daele for his help and advice when visiting the Department of Mathematics of the University of Leuven in Belgium during several times the past years. Shuanhong Wang  was partially supported  by the National Natural Science Foundation of China (Grant No. 11871144).
\nl

\section{\hspace{-17pt}. Preliminaries}  \label{s:prel} 

As we explained in the introduction, the most important parts of this note are  Sections \ref{s:dual1} and \ref{s:dual2}. But it would be impossible to explain the main ideas without recalling first some of the notions and results we plan to evaluate there. This is done in the present section, where we recall the notions of a regular weak multiplier Hopf algebra and a regular multiplier Hopf algebroid, together with their intimate relation. For this we essentially follow the content of \cite{T-VD2}. We will not give details as they can be found there.
\ssnl
In the next section we consider integrals and see how they give rise to integrals on the associated multiplier Hopf algebroid. To do so, it is also important to have the basic notions and properties as collected in this section.
\ssnl
We first recall briefly the concept of a (regular) weak multiplier Hopf algebra.
\ssnl
Let $A$ be a non-degenerate and idempotent algebra. Let $\Delta:A\to M(A\ot A)$ be a regular and full coproduct on this algebra. The  {\it canonical maps} $T_1$, $T_2$, $T_3$ and $T_4$ are defined on $A\ot A$ by
\begin{align}
T_1(a\ot b)&=\Delta(a)(1\ot b)
\qquad\quad\text{and}\quad\qquad
T_2(a\ot b)=(a\ot 1)\Delta(b) \label{eqn:1.1}\\
T_3(a\ot b)&=(1\ot b)\Delta(a)
\qquad\quad\text{and}\quad\qquad
T_4(a\ot b)=\Delta(b)(a\ot 1)\label{eqn:1.2}
\end{align}
when $a,b\in A$.  Regularity of the coproduct means that all these maps have range in $A\ot A$. These canonical maps play an important role further in this text.
\ssnl
It is assumed that the coproduct admits a \emph{counit}. This condition is related with the fullness property. Fullness of the coproduct means that the legs of $\Delta$ are all of $A$. This does not imply automatically that a counit exists. On the other hand, uniqueness of the counit (if it exists) does follow from the fullness of the coproduct.
\ssnl
One of the main features of a weak multiplier Hopf algebra is formulated in the following result.

\prop
There is a unique idempotent $E\in M(A\ot A)$ such that
\begin{equation*}
 E(A\ot A)=\Delta(A)(A\ot A)
\qquad\quad\text{and}\qquad\quad
(A\ot A)E=(A\ot A)\Delta(A).
\end{equation*}
It is the smallest idempotent element in $M(A\ot A)$ such that $\Delta(a)=E\Delta(a)=\Delta(a)E$ for all $a\in A$. It satisfies $(\Delta\ot\iota)E=(\iota\ot\Delta)E$, as well as
\begin{equation}
(\Delta\ot\iota)E=(E\ot 1)(1\ot E)=(1\ot E)(E\ot 1).\label{eqn:1.3}
\end{equation}
Recall that we use $\iota$ for the identity map.
\eprop
To give a meaning to these formulas, one has to extend the involved maps to the multiplier algebras. But this can be done precisely because of the first property in this proposition.
\ssnl
Another important result in the theory is the existence of an {\it antipode}.

\prop
There is a bijective linear map $S:A\to A$ such that the maps $R_1$ and $R_2$, defined on $A\ot A$ by
\begin{align*}
R_1(a\ot b)&=\sum_{(a)}a_{(1)}\ot S(a_{(2)})b \\
R_2(a\ot b)&=\sum_{(b)}aS(b_{(1)})\ot b_{(2)},
\end{align*}
are well-defined maps from $A\ot A$ to itself and so that $R_1$ and $R_2$ are generalized inverses of the canonical maps $T_1$ and $T_2$ respectively.
\eprop

We are using the Sweedler notation here. This is not obvious, but has been discussed in the original papers.
\ssnl
That $R_1$ is a generalized inverse of $T_1$ means that
\begin{equation*}
T_1R_1T_1=T_1
\qquad\quad\text\quad\qquad
R_1T_1R_1=R_1
\end{equation*}
and similarly for $T_2$ and $R_2$.
\ssnl
In general, such generalized inverses are not unique. They are determined by a choice of the projection map on the kernel and one on the range of the corresponding maps. Such a choice is part of the axioms of a weak multiplier Hopf algebra and it depends only on the canonical multiplier $E$. It follows that also the map $S$, as in the above proposition, is uniquely determined by these properties. It is called the antipode of the weak multiplier Hopf algebra $(A,\Delta)$. As expected it is an anti-isomorphism of the algebra $A$ and it also flips the coproduct $\Delta$.
\ssnl
It is understood that the above result is only true (as it is formulated) for a \emph{regular} weak multiplier Hopf algebra. Recall that a weak multiplier Hopf algebra is said to be regular if $(A,\Delta^{\text{cop}})$ (or equivalently $(A^{\text{op}},\Delta)$) also satisfies the axioms of a weak multiplier Hopf algebra. If we replace $\Delta$ by $\Delta^{\text{cop}}$, or $A$ by $A^{\text{op}}$, the antipode $S$ is replaced by the inverse $S^{-1}$. Then it is easy to find formulas for the appropriate generalized inverses $R_3$ and $R_4$ of $T_3$ and $T_4$ respectively.
\ssnl
For the ranges of the canonical maps we have the following formulas:

\prop
\begin{align*} T_1R_1(a\ot b)&=E(a\ot b)
\quad\qquad\text{and}\quad\qquad
T_2R_2(a\ot b)=(a\ot b)E\\
T_3R_3(a\ot b)&=(a\ot b)E
\quad\qquad\text{and}\quad\qquad
T_4R_4(a\ot b)=E(a\ot b)
\end{align*}
for all $a,b\in A$.
\eprop

We also have formulas for the kernels of the canonical maps.

\prop\label{prop:1.4}
There exists idempotents  $F_1, F_2, F_3$ and $F_4$ in $M(A\ot A^{\text{op}})$ given by
\begin{align*}
 F_1&=(\iota\ot S)E  \qquad\quad\text{and}\qquad\quad F_3=(\iota\ot S^{-1})E \\
	F_2&=(S\ot \iota)E \qquad\quad\text{and}\qquad\quad F_4=(S^{-1}\ot \iota)E.
\end{align*}
They give the kernels of the four canonical maps in the sense that
\begin{align*}
 R_1T_1(a\ot b)&=(a\ot 1)F_1(1\ot b)
\quad\quad\text{and}\quad\quad
R_2T_2(a\ot b)=(a\ot 1)F_2(1\ot b)\\
R_3T_3(a\ot b)&=(1\ot b)F_3(a\ot 1)
\quad\quad\text{and}\quad\quad
R_4T_4(a\ot b)=(1\ot b)F_4(a\ot 1)
\end{align*}
for all $a,b$ in $A$.
\eprop

\bf The source and target maps \rm
\nl
Having formulated the basic results about regular weak multiplier Hopf algebras, we now recall the {\it source} and {\it target maps}, the {\it source and target algebras} and their main properties needed further in this paper.

\defin\label{defin:1.5}
The source and target maps $\varepsilon_s$ and $\varepsilon_t$ from $A$ to $M(A)$ are defined by
\begin{equation*}
\varepsilon_s(a)=\sum_{(a)}S(a_{(1)})a_{(2)}
\qquad\quad\text{and}\quad\qquad
\varepsilon_t(a)=\sum_{(a)}a_{(1)}S(a_{(2)}).
\end{equation*}
We denote $\varepsilon_s(A)$ by $B$ and $\varepsilon_t(A)$ by $C$. These algebras are called the \emph{source} and \emph{target algebras}.
\edefin

These algebras satisfy the following properties.

\prop\label{prop:1.6}
We have $E\in M(B\ot C)$ and $E$ is a separability idempotent in the sense of \cite{VD4.v1}. The associated antipodal maps  $S_B:B \to C$ and $S_C:C \to B$ are precisely the restriction of the antipode (or rather of its extension to the multiplier algebra $M(A)$). This means that
\begin{equation*}
E(x\ot 1)=E(1\ot S(x))
\qquad\quad\text{and}\quad\qquad
(1\ot y)E=(S(y)\ot 1)E
\end{equation*}
for all $x\in B$ and $y\in C$.
\eprop

We make a comment on the notations.

\opm\label{opm:1.7}
i) In the papers on separability idempotents \cite{VD4.v1} and \cite{VD4.v2} we denoted the elements of the underlying algebras $B$ and $C$ with $b,b', \ldots$ and $c,c', \ldots$. However, in the context of a weak multiplier Hopf algebra, we use symbols $a,b,c$ for elements in the algebra $A$. Therefore, in this note, like in earlier papers on weak multiplier Hopf algebras, we  generally use $x$ for elements in $B$ and $y$ for elements in $C$. Remark however, in \cite{VD-W5}, it was the other way around. The convention was changed in the more recent paper \cite{T-VD2} because we felt that this is more natural.
\ssnl
ii) In what follows we will systematically use the symbols $S_B$ and $S_C$ for the anti-isomorphisms obtained by the restriction to $B$ and $C$ resp. of the extension of the antipode from $A$ to $M(A)$.
\eopm

Also the following is important.

\opm\label{opm:1.8}
Observe that $E(1\ot y)\in B\ot C$ for all $y\in C$ and if we use a Sweedler type notation $E_{1}\ot E_{2}$ for $E$ we have $S_B(E_{1})E_{2}y=y$ when $y\in C$. Similarly, $(x\ot 1)E\in B\ot C$ and $xE_{1}S_C(E_{2})=x$ for all $x\in B$. We write these properties as
\begin{equation*}
S_B(E_{1})E_{2}=1
\qquad\quad\text{and}\qquad\quad
E_1S_C(E_2)=1.
\end{equation*}
See Proposition 1.8 in \cite{VD4.v2}.
\eopm
\ssnl
We not only have that $E$ sits in $M(B\ot C)$  but moreover that $B$ is the left leg and $C$ the right leg of $E$ (in an appropriate sense). Then it follows from Equation (\ref{eqn:1.3}) that these algebras commute with each other and that in fact
\begin{equation}
\Delta(x)=E(1\ot x)=(1\ot x)E
\tussenen
\Delta(y)=(y\ot 1)E=E(y\ot 1)\label{eqn:1.4a}
\end{equation}
for all $x\in B$ and $y\in C$.
\ssnl
We will also need the linear functionals $\mu_B$ and $\mu_C$ on $B$ and $C$ respectively as we recall in the following proposition.

\prop\label{prop:1.8}
There exist unique linear functionals $\mu_B$ and $\mu_C$ on $B$ and $C$ satisfying
\begin{equation}
(\mu_B\ot\iota)E=1
\tussenen
(\iota\ot \mu_C)E=1.\label{eqn:1.5a}
\end{equation}
We have $\mu_C\circ S_B=\mu_B$ and $\mu_B\circ S_C=\mu_C$. Finally, they satisfy the following KMS properties. For all $x,x'\in B$ and $y,y'\in C$ we have
\begin{equation}
\mu_B(x'x)=\mu_B(S_CS_B(x)x')
\tussenen
\mu_C(yy')=\mu_C(y'S_BS_C(y))\label{eqn:1.6a}
\end{equation}
\eprop

The Equations (\ref{eqn:1.5a}) are valid in the multiplier algebras of $C$ and $B$ respectively. The formulas in (\ref{eqn:1.6a}) tell us that the modular automorphism $\sigma_B$ of $\mu_B$ is the inverse of the composition of the antipodal maps while the modular automorphism $\sigma_C$ is the composition of the antipodal maps.
\snl
Finally, we mention that the algebras $B$ and $C$ have local units (see Proposition 1.10 in \cite{VD4.v2}).
The same is true for the original algebra $A$ (see Proposition 4.9 in \cite{VD-W4}).
\ssnl
For details we refer to the original papers \cite{VD-W3, VD-W4, VD-W5} as well as to the preliminary section of \cite{T-VD2} for an expanded version of the notions and properties above.
\nl
\bf The associated multiplier Hopf algebroid \rm
\nl
A regular multiplier Hopf algebroid is a difficult and complex concept. Because here we are interested in the multiplier Hopf algebroid that comes from a weak multiplier Hopf algebra, we can avoid going deeply into the theory. This will also not be necessary for the forthcoming discussions in this note.
\ssnl
The easiest and most transparent way to pass to the multiplier Hopf algebroid in this case is to describe the new canonical maps and how they are constructed from the original ones, given in Equations (\ref{eqn:1.1}) and (\ref{eqn:1.2}). These new maps are now bijections between balanced tensor products.
\ssnl
To define these balanced tensor products, we first introduce the following module structures. See e.g.\ Item 2.1 in \cite{T2}.

\notat\label{notat:1.8}
In the first place, we define the anti-isomorphism $t_B:B\to C$ and $t_C:C\to B$ simply as the inverses of the maps $S_C:C\to B$ and $S_B:B\to C$ respectively. Recall that these maps are the restrictions of the antipode $S$ after it has been extended to the multiplier algebra $M(A)$ (see Remark \ref{opm:1.7}.ii).
\enotat

\notat
i) We make $A$ into a left $B$-module in two ways. First we have $(x,a)\mapsto xa$ when $a\in A$ and $x\in B$. We use $_BA$ for this module. Secondly we have $(x,a)\mapsto at_B(x)$ for $a\in A$ and $x\in B$. We use ${^B\hskip -3pt}A$ for this module.
\ssnl
ii) Further we make $A$ into a right $B$-module in two ways. We use $A_B$ for the action $(a,x)\mapsto ax$ while we use $A^B$ for the action $(a,x)\mapsto t_B(x)a$.
\ssnl
iii) In a completely similar way, we have two left and two right $C$-modules. For the left modules we use $_CA$ and $^C\hskip -3ptA$ while for the right modules we use $A_C$ and $A^C$. Of course, now the map $t_C$ is used instead of $t_B$ to define the modules $^C\hskip -3ptA$ and $A^C$.
\enotat

These are the notations used in \cite{T-VD1}. In this note, we will in most cases use the antipodal maps $S_C$ and $S_B$ and not their inverses $t_B$ and $t_C$.
\ssnl
With these module structures come exactly eight balanced tensor products, but we will only need six of them. The first two are
\begin{equation*}
A_B\ot {_BA}
\tussenen
A_C\ot {_CA}.
\end{equation*}
So in $A_B\ot {_BA}$ we have $ax\ot b=a\ot xb$ for all $a,b\in A$ and $x\in B$. Similarly in $A_C\ot {_CA}$ it holds $ay\ot b=a\ot yb$ for all $a,b\in A$ and $y\in C$. Remark that in \cite{T-VD2} we have used $A\ot_s A$ for the first one and $A\ot_t A$ for the second one.
\ssnl
The next two are derived from these by applying the flip map. They are
\begin{equation*}
A^B\ot {^B\hskip-3ptA}
\tussenen
A^C\ot {^C\hskip -3pt A}.
\end{equation*} Now in $A^B\ot {^B\hskip-3ptA}$ we have $t_B(x)a\ot b=a\ot bt_B(x)$ for all $a,b\in A$ and $x\in B$. Equivalently $ya\ot b=a\ot by$ for all $a,b\in A$ and $y\in C$. This indeed is obtained from $A_C\ot {_CA}$ by applying the flip. Similarly for $A^C\ot {^C\hskip-3ptA}$ that is obtained from $A_B\ot {_BA}$ by the flip map. In \cite{T-VD2} these spaces are denoted with $A\ot^s A$ and $A\ot^t A$ respectively.
\ssnl
Finally, the last two are of a different nature:
\begin{equation*}
A_B\ot {^B\hskip -3pt A}
\tussenen
A^C\ot {_CA}.
\end{equation*}
In $A_B\ot {^B\hskip-3ptA}$ we have $ax\ot b=a\ot bt_B(x)$ while in $A^C\ot {_CA}$ we have $t_C(y)a\ot b=a\ot yb$ for all $a,b\in A$ and $x\in B$ and $y\in C$. Compared with the notations in \cite{T-VD2} we have there $A\ot_r A$ for $A_B\ot {^B\hskip-3ptA}$ and $A\ot_\ell A$ for $A^C\ot {_CA}$. Observe that $A\ot_r A$ is defined in \cite{T-VD2} so that $a\ot by=aS_C(y)\ot b$ for all $y\in C$ which is the same as $ax\ot b=a\ot bt_B(x)$
for all $x\in B$, precisely because $t_B=S_C\inv$.
\ssnl
The remaining two are not used in this theory.
\ssnl
In what follows we will use $\pi$ for all the projection maps from $A\ot A$ to the various balanced tensor products.

\prop\label{prop:1.10}
There are bijective maps,
\begin{align}
&\mathfrak T_1:A_B\ot {_BA}\to A^C\ot {_CA}
&\text{and}\quad\quad
&\mathfrak T_2: A_C\ot {_CA}\to A_B\ot {^B\hskip -3pt A} \label{eqn:1.4}\\
&\mathfrak T_3: A^C \ot{^C\hskip -3pt A}\to  A_B\ot {^B\hskip -3pt A}
&\text{and}\quad\quad
&\mathfrak T_4:A^B\ot {^B\hskip -3pt A}\to A^C\ot {_CA} \label{eqn:1.5}
\end{align}
satisfying
\begin{align}
&\mathfrak T_1\circ \pi=\pi\circ T_1
\qquad\qquad\text{and}\qquad\qquad
\mathfrak T_2\circ \pi=\pi\circ T_2\\
&\mathfrak T_3\circ\pi=\pi\circ T_3
\qquad\qquad\text{and}\qquad\qquad
\mathfrak T_4\circ\pi=\pi\circ T_4.
\end{align}
We have the obvious quotient maps in each of these formulas.
\eprop

In other words, these maps give commutative diagrams with $T_1$ and $T_2$ for maps in the Equations (\ref{eqn:1.4}) respectively, and with $T_3$ and $T_4$ for the maps in the Equations (\ref{eqn:1.5}).
\ssnl
The bijectivity of the four maps $(\mathfrak T_i)$ is a consequence of the properties of the range and kernel of the original canonical maps $(T_i)$. This has been shown in detail in Section 2 of \cite{T-VD2}.
\ssnl
These maps are the canonical maps of the regular multiplier algebroid associated with the original weak multiplier Hopf algebra. They contain all information necessary for the further discussion in this note. More can be found in the first three sections of \cite{T-VD2}, in particular when we are  using Theorem 3.5 of that paper.

\notat
We will use $\mathfrak A$ to denote the regular multiplier Hopf algebroid associated with the given regular weak multiplier Hopf algebra $A$. The maps $\mathfrak T_1$, $\mathfrak T_2$, $\mathfrak T_3$ and $\mathfrak T_4$ are the canonical maps of $\mathfrak A$.
\enotat

The coproducts $\Delta_B$ and $\Delta_C$ of $\mathfrak A$ are maps defined on $A$. The target is not so easy to describe and we refer to \cite{T-VD1} for a precise definition.  The coproducts are completely determined (and in fact implicitly defined) by these canonical maps  and we have
\begin{align}
&\mathfrak T_1(a\ot b)=\Delta_C(a)(1\ot b)
\qquad\quad\text{and}\qquad\quad
\mathfrak T_2(a\ot b)=(a\ot 1)\Delta_B(b)\label{eqn:1.9}\\
&\mathfrak T_3(a\ot b)=(1\ot b)\Delta_B(a)
\qquad\quad\text{and}\qquad\quad
\mathfrak T_4(a\ot b)=\Delta_C(b)(a\ot 1).\label{eqn:1.10}
\end{align}
Further in the paper, we always work with these canonical maps. The above statement only serves a motivational purpose.
\nl
Let us finish this section with a \emph{dictionary}. A regular multiplier Hopf algebroid, as it appears in the papers \cite{T-VD1}, \cite{T1} and \cite{T2} is in fact a tuple $(A,B,C,t_B,t_C,\Delta_B,\Delta_C)$. For the multiplier Hopf algebroid $\mathfrak A$ associated to the weak multiplier Hopf algebra $(A,\Delta)$ as described above, $A$, $B$ and $C$ of the tuple are precisely the original algebras $A$, $B$ and $C$ as in these preliminaries. The maps $t_B$ and $t_C$ are the inverses of the antipode as we explained in Notation \ref{notat:1.8}. We get
\begin{equation*}
t_B(x)=S_C\inv(x)
\tussenen
t_C(y)=S_B\inv(y)
\end{equation*}
for $x\in B$ and $y\in C$. Remark again that we use $S_B$ and $S_C$ for the restrictions to $B$ and $C$ respectively of the extension to the multiplier algebra of the map $S$.
\ssnl
Finally remark that in the original papers on multiplier Hopf algebroids, these canonical maps are denoted differently. One has
\begin{align*}
&\mathfrak T_1=T_\rho
\tussenen
\mathfrak T_2={_\lambda}T\\
&\mathfrak T_3={_\rho}T
\tussenen
\mathfrak T_4=T_\lambda.
\end{align*}

We plan to avoid these notations for the canonical maps by using instead the right hand sides in the formulas (\ref{eqn:1.9}) and (\ref{eqn:1.10}).
\ssnl
The counital maps are denoted by $\varepsilon_B$ and $_C\varepsilon$ in the original papers on multiplier Hopf algebroids. For the multiplier algebroid $\mathfrak A$ obtained from the weak multiplier Hopf algebra, we precisely have
\begin{equation*}
\varepsilon_B=\varepsilon_s
\tussenen
 _C\varepsilon=\varepsilon_t
\end{equation*}
as expected. Finally, the antipode of the multiplier algebroid $\mathfrak A$ is nothing else but the original antipode.

\section{\hspace{-17pt}. From algebraic quantum groupoids to multiplier Hopf algebroids with integrals}  \label{s:from}

In the previous section, we have seen how a regular weak multiplier Hopf algebra gives rise to a regular multiplier Hopf algebroid. This has been shown in detail in \cite{T-VD2} and we just have collected the main steps from that paper needed further for this note.
\ssnl
In this section, we start again with a regular weak multiplier Hopf algebra $(A,\Delta)$. Now we assume that it has \emph{a faithful set of integrals}. As mentioned in the introduction, we then call it an algebraic quantum groupoid. The aim of this section is to show how integrals on $(A,\Delta)$ give rise to integrals on the associated algebroid. This has been shown in Theorem 5.1.4 of Section 5.1 of \cite{T2} (in the case of a single faithful integral). Our approach here is slightly different and we will comment on the connection between this approach and the one in \cite{T2}.
\nl
First recall Definition 1.1  from \cite{VD-W6}. Here again $B$ and $C$ denote the source and the target algebras of $(A,\Delta)$ as in Definition \ref{defin:1.5} of the previous section.

\defin\label{defin:2.1}
A linear functional $\varphi$ on $A$ is called  \emph{left invariant} if $(\iota\ot\varphi)\Delta(a)\in M(C)$ for all $a\in A$. A \emph{non-zero} left invariant functional is called a \emph{left integral}. Similarly, a linear functional $\psi$ on $A$  is called \emph{right invariant}  if $(\psi\ot\iota)\Delta(a)\in M(B)$ for all $a\in A$ and a non-zero right invariant functional is called a \emph{right integral}.
\edefin

\opm In \cite{VD-W5} we have treated the source and target algebras in the more general, not necessarily regular case. Observe however, that for regular weak multiplier Hopf algebras, the algebras $A_t$ and $A_s$ as defined in Notation 2.5 of \cite{VD-W5} coincide with the multiplier algebras $M(C)$ and $M(B)$ respectively, see Proposition 2.16 of \cite{VD-W5}.
\eopm

We recall the following fundamental property of the  integrals.

\prop\label{prop:2.3}
Assume that $\varphi$ is a left invariant functional and that $\psi$ is a right invariant functional on $A$. Then we have
\begin{align}
(\iota\ot\varphi)\Delta(a) &=(\iota\ot\varphi)(F_2(1\ot a))=(\iota\ot\varphi)((1\ot a)F_4) \label{eqn:2.1}\\
(\psi\ot\iota)\Delta(a)&=(\psi\ot\iota)((a\ot 1)F_1)=(\psi\ot\iota)(F_3(a\ot 1)).\label{eqn:2.2}
\end{align}
for all $a\in A$.
\eprop

The elements $F_i$ are defined in Proposition \ref{prop:1.4}. This result is found in Proposition 1.4 of \cite{VD-W6} where it is in fact a reformulation of Proposition 1.3 in \cite{VD-W6}, at least partly. The proof is also found in Proposition 3.7 of \cite{K-VD}.
\ssnl
It follows from these formulas that in fact $(\iota\ot\varphi)\Delta(a)$ belongs to $C$ for all $a$ and not only to the multiplier algebra $M(C)$ as originally assumed. Similarly $(\psi\ot\iota)\Delta(a)\in B$ for all $a$. Indeed, from the properties of $E$ (see e.g.\ Remark \ref{opm:1.8}) and the formula $F_1=(\iota\ot S)E$, we have that $(a\ot 1)F_1\in A\ot B$ so that  $(\psi\ot\iota)((a\ot 1)F_1)\in B$. Similarly for the three other cases.
\ssnl
In terms of the Sweedler notation, these  formulas read as
\begin{align}
(\iota\ot\varphi)\Delta(a) &=\sum_{(a)}a_{(1)}S(a_{(2)})\varphi(a_{(3)})=\sum_{(a)}S\inv(a_{(2)})a_{(1)}\varphi(a_{(3)})\label{eqn:2.3}\\
(\psi\ot\iota)\Delta(a)&=\sum_{(a)}\psi(a_{(1)})S(a_{(2)})a_{(3)}=\sum_{(a)}\psi(a_{(1)})a_{(3)}S\inv(a_{(2)})\label{eqn:2.4}
\end{align}
\opm
i) In Definition 5.1.3 of \cite{T2} a linear functional $\varphi$ on a regular weak multiplier Hopf algebra is called a left integral if it satisfies the first equality of Equation (\ref{eqn:2.1}) for all $a$ in $A$. Similarly in \cite{T2}, a linear functional $\psi$ is called a right integral if it satisfies the first equality of Equation (\ref{eqn:2.2}) for all $a$.
\ssnl
ii) As it is clear that these properties imply that they are left, resp.  right invariant in the sense of Definition \ref{defin:2.1}, we see that the two definitions yield the same objects. Just observe that we only call an invariant functional an integral when it is non-zero.
\ssnl
iii) We prefer our definition over the one in \cite{T2} because it is simpler and more direct. It can be formulated just in terms of the coproduct. Strictly speaking, one even does not need the source and target maps for defining the algebras $M(B)$ and $M(C)$. They are characterized as the algebras of elements $x, y\in M(A)$ resp. satisfying
\begin{equation*}
\Delta(ax)=\Delta(a)(1\ot x)
\tussenen
\Delta(ya)=(y\ot 1)\Delta(a)
\end{equation*}
for all $a$. It should be mentioned however that in the theory of weak Hopf algebras, the Equations (\ref{eqn:2.3}) and (\ref{eqn:2.4}) are used to define invariance. In the case of weak multiplier Hopf algebras, this would not be an obvious approach. This is discussed in Section 1 of \cite{VD-W6} where we introduce the notion of integrals for regular weak multiplier Hopf algebras.
\eopm

In Definition 5.1.3 of  \cite{T2} a regular weak multiplier Hopf algebra is called \emph{a regular weak multiplier Hopf algebra with integrals} only when it has a left integral $\varphi$  such that
each of the  sets of elements
\begin{equation}
(\varphi\ot\iota)((a\ot 1)E) \tussenen (\varphi\ot\iota)(E(a\ot 1))\label{eqn:2.5}
\end{equation}
where $a\in A$, is equal to $C$,  and if it has is a right integral $\psi$ such that the sets
\begin{equation}
(\iota\ot\psi)((1\ot a)E) \tussenen (\iota\ot\psi)(E(1\ot a))\label{eqn:2.6}
\end{equation}
where $a\in A$ both are all of $B$.

\opm
i) The terminology is a bit strange. With this terminology, it is not sufficient for a regular weak multiplier Hopf algebra to have integrals to be called a regular multiplier Hopf algebra with integrals.
\ssnl
ii) As a matter of fact, it can be shown that these conditions on the integrals just mean that they are faithful. This has been shown using standard techniques in Lemma 1.15, Lemma 1.16 and Proposition 1.17 of \cite{VD-W6}. One direction is easy and holds in fact for any linear functional. Indeed, as soon as $\varphi$ is a faithful linear functional on $A$ we will have that all of $C$ is spanned by elements as above in Equation (\ref{eqn:2.5}). Similarly for $\psi$. This is shown in Lemma 1.15 of \cite{VD-W6}.
\ssnl
For the other direction, it is needed that the functionals are integrals. It is shown in Lemma 1.16 of \cite{VD-W6}. The proof is inspired by the arguments that are used to show that an integral on a regular multiplier Hopf algebra is automatically faithful.
\ssnl
The combination of the two results is then formulated in Proposition 1.17 of \cite{VD-W6}. So in the terminology of \cite{T2},  a regular weak multiplier Hopf algebra with integrals is in fact a regular weak multiplier Hopf algebra with a single faithful integral. Remark further that if there is a faithful left integral, there is also a faithful right integral. In particular, it means that one of the conditions above will be sufficient as it will imply the other one.
\eopm

In \cite{VD-W6} we do not assume the existence of a single faithful integral however to construct the dual, we only need a \emph{faithful set} of left integrals. This is indeed a weaker condition because it can happen that there is such a faithful set but not a single faithful integral, see e.g.\ \cite{I-K}. We are indebted to G.\ B\"ohm for drawing our attention to this example.
\ssnl
Because the theory in \cite{T2} is only developed in the case of a single faithful integral, we will stick to that case when refering to \cite{T2}. We believe however that this theory can be generalized to the case of a faithful set of integrals as in the duality theory of weak multiplier Hopf algebras with integrals, developed in \cite{VD-W6}. See also the discussion on this topic in Section \ref{s:concl}.

\nl
\bf The passage to the multiplier Hopf algebroid \rm
\nl
In what follows, we assume that $(A,\Delta)$ is a weak multiplier Hopf algebra and that $\varphi$  is a left integral and $\psi$ a right integral on $A$ to begin with. For the moment, we do not yet  assume that there is a faithful set of integrals. We will do that later.
\snl
Consider the multiplier Hopf algebroid $\mathfrak A$ associated to $(A,\Delta)$ as in the previous section.

\defin\label{defin:2.5}
Define the linear maps $\Phi:A\to C$ and $\Psi:A\to B$ by
\begin{equation*}
\Phi(a)=(\iota\ot\varphi)\Delta(a)
\tussenen
\Psi(a)=(\psi\ot\iota)\Delta(a)
\end{equation*}
for $a\in A$.
\edefin

We will now show that $\Phi$ is a partial left integral and that  $\Psi$ is a partial right integral in the sense of Definition 2.2.1 of \cite{T2}.

\prop\label{prop:2.7a}
The linear map $\Phi$ is a module map from $_CA_C$ to $_CC_C$ in the sense that
\begin{equation*}
\Phi(y_1ay_2)=y_1\Phi(a)y_2
\end{equation*}
for all $a$ in $A$ and $y_1,y_2\in C$. Similarly the map $\Psi$ is a module map from $_BA_B$ to $_BB_B$.
\eprop
\bew
This follows immediately from the the Equations  (\ref{eqn:1.4a}) in the previous section. Indeed, we have e.g.\ that $\Delta(ya)=(y\ot 1)\Delta(a)$ and hence $\Phi(ya)=y\Phi(a)$ for all $a\in A$ and $y\in C$. Similar arguments are used for the three other cases.
\ebew

\prop\label{prop:2.6}
We have
\begin{align}
&m_1(\Psi\ot\iota)(\Delta_C(a)(1\ot b))=\Psi(a)b\label{eqn:2.7}\\
&m_2(\iota\ot \Phi)((a\ot 1)\Delta_B(b))=a\Phi(b)\label{eqn:2.8}
\end{align}
for all $a,b\in A$. For these two equations, we use the (well-defined) maps $m_1:a\ot b\mapsto S(a)b$ from $A^C\ot {_CA}$ to $A$ and $m_2:a\ot b\mapsto aS(b)$ from $A_B\ot{^B\hskip -3pt A}$ to $A$.
\eprop

\bew
Using the Sweedler notation we find using Equation (\ref{eqn:2.4})
\begin{equation*}
(\Psi\ot\iota)(\Delta_C(a)(1\ot b))=\sum_{(a)}\Psi(a_{(1)})\ot a_{(2)}b=\sum_{(a)}\psi(a_{(1)})a_{(3)}S\inv(a_{(2)})\ot a_{(4)}b
\end{equation*}
and if we apply the map $m_1$ we get, using Equation (\ref{eqn:2.4})
\begin{align*}
m_1((\Psi\ot\iota)(\Delta_C(a)(1\ot b)))
&=\sum_{(a)}\psi(a_{(1)})S(a_{(3)}S\inv(a_{(2)})) a_{(4)}b\\
&=\sum_{(a)}\psi(a_{(1)})a_{(2)}S(a_{(3)})a_{(4)}b\\
&=\sum_{(a)}\psi(a_{(1)})a_{(2)}b=\Psi(a)b.
\end{align*}
This proves Equation (\ref{eqn:2.7}). Similarly we have, using again Equation \ref{eqn:2.3},
\begin{equation*}
(\iota\ot \Phi)((a\ot 1)\Delta_B(b))=\sum_{(b)}ab_{(1)}\ot \Phi(b_{(2)})=\sum_{(b)}ab_{(1)}\ot S\inv(b_{(3)})b_{(2)}\varphi(b_{(4)})
\end{equation*}
and if we apply the map $m_2$ we find
\begin{align*}
m_2(\iota\ot \Phi)((b\ot 1)\Delta_B(a))
&=\sum_{(b)}ab_{(1)}S(S\inv(b_{(3)})b_{(2)})\varphi(b_{(4)})\\
&=\sum_{(b)}ab_{(1)}S(b_{(2)})b_{(3)}\varphi(b_{(4)})\\
&=\sum_{(b)}ab_{(1)}\varphi(b_{(2)})=a\Phi(b).
\end{align*}
This proves Equation (\ref{eqn:2.8}).
\ebew

The proof can be formulated without the use of the Sweedler notation by using the formulas from Proposition \ref{prop:2.3}. The crucial one for the first case is the formula
\begin{equation*}
(\psi\ot\iota)\Delta(a)=( \psi\ot\iota)(F_3(a\ot 1))
\end{equation*}
and for the second case, it is
\begin{equation*}
(\iota\ot\varphi)\Delta(a)=(\iota\ot\varphi)((1\ot a)F_4).
\end{equation*}
We have a similar result for the two other maps.

\prop\label{prop:2.7}
We have
\begin{align}
&m_3(\Psi\ot\iota)((1\ot b)\Delta_B(a))=b\Psi(a)\label{eqn:2.9}\\
&m_4(\iota\ot \Phi)(\Delta_C(b)(a\ot 1))=\Phi(b)a\label{eqn:2.10}
\end{align}
for all $a,b\in A$. For these two equations, we use the (well-defined) maps $m_3:a\ot b\mapsto bS\inv(a)$ from $A_B\ot {^B\hskip -3pt A}$ to $A$ and $m_4:a\ot b\mapsto S\inv(b)a$ from $A^C\ot{_CA}$ to $A$.
\eprop
\bew
In the first case we have, using Equation (\ref{eqn:2.4}),
\begin{equation*}
(\Psi\ot\iota)((1\ot b)\Delta_B(a))=\sum_{(a)}\Psi(a_{(1)})\ot ba_{(2)}=\sum_{(a)}\psi(a_{(1)})S(a_{(2)})a_{(3)}\ot ba_{(4)}
\end{equation*}
and if we apply the map $m_3$ we find
\begin{align*}
m_3(\Psi\ot\iota)((1\ot b))\Delta_B(a))
&=\sum_{(a)}\psi(a_{(1)}) ba_{(4)}S\inv(a_{(3)})a_{(2)}\\
&=\sum_{(a)}\psi(a_{(1)}) ba_{(2)}=b\Psi(a).
\end{align*}
For the second formula we use Equation (\ref{eqn:2.3}) again and we find
\begin{equation*}
(\iota\ot \Phi)(\Delta_C(b)(a\ot 1))=\sum_{(b)}b_{(1)}a\ot \Phi(b_{(2)})=\sum_{(b)}b_{(1)}a\ot b_{(2)}S(b_{(3)}) \varphi(b_{(4)}).
\end{equation*}
When we apply $m_4$ we find
\begin{align*}
m_4((\iota\ot \Phi)(\Delta_C(b)(a\ot 1)))
&=\sum_{(b)}S\inv(b_{(2)}S(b_{(3)})) \varphi(b_{(4)})b_{(1)}a\\
&=\sum_{(b)}  \varphi(b_{(4)})b_{(3)}S\inv(b_{(2)}b_{(1)}a=\Phi(b)a.
\end{align*}
\vskip -20pt
\ebew

Also here, the results can be obtained without the use of the Sweedler notation, but using the original formulas in Proposition \ref{prop:2.3}.
\ssnl
We include some comments and compare these formulas with the formulas (LI1), (LI2) and (RI1), (RI2) of Section 2.2 in \cite{T2}.

\opm\label{opm:2.10}
i) First consider the expressions
\begin{equation*}
(\Psi\ot\iota)((1\ot b)\Delta_B(a))
\tussenen
(\iota\ot\Phi)(\Delta_C(b)(a\ot 1))
\end{equation*}
that we consider in Proposition \ref{prop:2.7}.
For the first one, we get an element in $B\ot A$ as sitting in $A_B\ot {^B\hskip -3pt A}$. We have the natural identification $x\ot b\mapsto x.b$ where $x.b=bS\inv(x)$. This identification is the restriction of the map $a\ot b\mapsto bS\inv(a)$. We denoted this map in Proposition \ref{prop:2.7} by $m_3$. From this point of view, we can rewrite Equation (\ref{eqn:2.7}) simply as
\begin{equation*}
(\Psi\ot\iota)((1\ot b)\Delta_B(a))=b\Psi(a).
\end{equation*}
This is the condition (RI2) from \cite{T2}.
\ssnl
For the second one we get an element in $A\ot C$, as sitting in $A^C\ot {_CA}$. Here we have a natural identification $a\ot y\mapsto a.y$ with $a.y=S\inv(y)a$. Then this identification is the restriction of the map $a\ot b\mapsto S\inv(b)a$. We denoted this map in Proposition \ref{prop:2.7a} by $m_4$. From this point of view, we can rewrite Equation (\ref{eqn:2.10}) simply as
\begin{equation*}
(\iota\ot\Phi)(\Delta_C(b)(a\ot 1)=\Phi(b)a.
\end{equation*}
This is condition (LI2) from \cite{T2}.
\ssnl
ii) The situation with the expressions
\begin{equation*}
(\Psi\ot\iota)(\Delta_C(a)(1\ot b))
\tussenen
(\iota\ot\Phi)((a\ot 1)\Delta_C(b))
\end{equation*}
that we consider in Proposition \ref{prop:2.6} is slightly different.
The first expression belongs to $B\ot A$ as sitting in  $A^C\ot {_CA}$. If we apply the antipode $S$ on the first leg, we find an element in $C\ot A$ as siting in $A_C\ot {_CA}$. Now we use the natural identification $y\ot a\mapsto ya$ and in combination with the antipode we precisely have the restriction of the $m_1$. So we can write
\begin{equation*}
(S_C\Psi\ot\iota)(\Delta_C(a)(1\ot b))=\Psi(a)b.
\end{equation*}
This gives us condition (RI1) of \cite{T2}.
\ssnl
Similarly  the second expression belongs to $A\ot C$ as sitting in $A_B\ot {^B\hskip -3pt A}$. We compose again with $S$ in the second leg and we get an expression in $A\ot B$, sitting in $A_B\ot {_BA}$. We have now the identification $a\ot x\mapsto ax$ and in combination with the antipode, we have the restriction of the map $m_2$. So we can write
\begin{equation*}
(\iota\ot S_B\Phi)((a\ot 1)\Delta_C(b))=a\Phi(b)
\end{equation*}
and we obtain condition (LI1) of \cite{T2}.
\eopm

The next step is to look for the other slice maps in item (4) of Definition 2.2.1 in \cite{T2}. They are obtained in the following proposition. We use the notations from \cite{T2}.

\prop\label{prop:2.11a}
 i) Define the maps $\Phi_B$ and $_B\Phi$ from $A$ to $B$ by
\begin{equation*}
\Phi_B(a)=(\iota\ot \varphi)((1\ot a)F_1)
\tussenen
_B\Phi(a)=(\iota\ot\varphi)(F_3(1\ot a)).
\end{equation*}
Then $\Phi_B(ax)=\Phi_B(a)x$ and $_B\Phi(xa)=x{_B\Phi}(a)$ for all $a$ in $A$ and $x$ in $B$. Moreover
$$\varphi(a)=\mu_B(\Phi_B(a))=\mu_B(_B\Phi(a))$$
for all $a\in A$. Here $\mu_B$ is the map satisfying $(\mu_B\ot\iota)E=1$ as in Proposition \ref{prop:1.8}.
\snl
ii) Define the maps $\Psi_C$ and $_C\Psi$ from $A$ to $C$ by
\begin{equation*}
\Psi_C(a)=(\psi\ot\iota)((a\ot 1)F_2)
\tussenen
_C\Psi(a)=(\psi\ot\iota)(F_4(a\ot 1)).
\end{equation*}
Then $\Psi_C(ay)=\Psi_C(a)y$ and $_C\Psi(ya)=y_C\Psi(a)$ for all $a$ in $A$ and $y$ in $C$. Moreover
$$\psi(a)=\mu_C(\Psi_C(a))=\mu_C(_C\Psi(a))$$ for all $a\in A$. Now $\mu_C$ satisfies $(\iota\ot\mu_C)E=1$ as in Proposition \ref{prop:1.8}.
\eprop

\bew
The proof is straightforward and uses the basic formulas from Proposition 1.6.
\ssnl
The formulas $\varphi(a)=\mu_B(\Phi_B(a))=\mu_B(_B\Phi(a))$ are a consequence of the first  formula in Equation (\ref{eqn:1.5a}) of Proposition \ref{prop:1.8} while, similarly the formulas $\psi(a)=\mu_C(\Psi_C(a))=\mu_C(_C\Psi(a))$ follow essentially from the second formula in Equation (\ref{eqn:1.5a}).
\ebew

If we combine all these results, we find the following.

\stel
Let $(A,\Delta)$ be a regular weak multiplier Hopf algebra with a faithful left integral $\varphi$ and a faithful right integral $\psi$. The maps $\Phi$ and $\Psi$ as defined in Definition \ref{defin:2.5} make the associated regular multiplier Hopf algebroid $\mathfrak A$ into a measured regular multiplier Hopf algebroid in the sense of Definition 2.2.1 of \cite{T2}.
\estel

\bew
i) In Proposition \ref{prop:2.7a} we have shown the $\Phi$ and $\Psi$ have the required module properties.
\ssnl
ii) We use the faithful functionals $\mu_B$ and $\mu_C$ from Proposition \ref{prop:1.8}. We have $\mu_C=\mu_B\circ S_C$ and $\mu_B=\mu_C\circ S_B$ taking care of the first part of Item (1) in Definition 2.2.1 of \cite{T2}. We have
\begin{equation*}
E(a\ot 1)=\sum_{(a)}a_{(1)}\ot a_{(2)}S(a_{(3)})
\end{equation*}
for all $a$. If we apply $\varepsilon\ot \mu_C$ we will find $\varepsilon(a)=\mu_C(\varepsilon_t(a))$. Similarly we get $\varepsilon(a)=\mu_B(\varepsilon_s(a))$ if we apply $\mu_B\ot\varepsilon$ on $(1\ot a)E$. This proves the second part of Item (1) in Definition 2.2.1 of \cite{T2} because $\varepsilon_s=\varepsilon_B$ and $\varepsilon_t={_C}\varepsilon$.
\ssnl
iii) In Proposition \ref{opm:2.10} we have seen that the conditions (LI) and (RI) of Definition 2.2.1 of \cite{T2} hold $\Phi$ and $\Psi$ respectively.
\ssnl
iv) The give left and right integrals $\varphi$ and $\psi$ satisfy $\varphi=\mu_C\circ\Phi$ and $\psi=\mu_B\circ \Psi$. Indeed, we have e.g.
\begin{equation*}
\Phi(a)=(\iota\ot\varphi)\Delta(a)=(\iota\ot\varphi((F_2(1\ot a))
\end{equation*}
where $F_2=(S_B\ot \iota)E$. If we apply $\mu_C$ and use that $\mu_C\circ S_B=\mu_B$ we find $\mu_C(\Phi(a))=\varphi(a)$. Similarly for the other formula. So also Item (3) of Definition 2.2.1 of \cite{T2} is satisfied.
\ssnl
v) Finally, we have the maps $_B\Phi, \Phi_B, {_C}\Psi, \Psi_C$ needed for Item (4) of Definition 2.2.1 of \cite{T2} with the required properties. This is shown in  Proposition \ref{prop:2.11a}. That these maps are surjective follows from the fact that the integrals $\varphi$ and $\psi$ are faithful and that $E$ is full. See the discussion in Remark \ref{opm:2.10}
\ebew

 In Section \ref{s:dual3}, where we treat special cases, we will come back to this construction. In the first item of that section, we start with a dual pair of regular weak multiplier Hopf algebra $A$ and $A'$ and consider the associated regular multiplier Hopf algebroids $\mathfrak A$ and $\mathfrak A'$. We treat the dual of a regular multiplier Hopf algebra with integrals in that context. We refer to Theorem \ref{stel:5.11} and Theorem \ref{stel:5.12}.

\section{\hspace{-17pt}. Pairing and duality of algebraic quantum groupoids}  \label{s:dual1}

Consider a regular weak multiplier Hopf algebra $(A,\Delta)$. If it has a faithful set of integrals we can consider the dual $(\widehat A,\widehat\Delta)$ as constructed in Section 2 of \cite{VD-W6}. It is again a regular weak multiplier Hopf algebra with integrals. The result is found in Theorem 2.15 and Theorem 2.19 of \cite{VD-W6}.
\ssnl
We can associate the regular multiplier Hopf algebroid $\mathfrak A$ to $A$ and we can do the same for the dual $\widehat A$. We will obtain a pairing of regular multiplier Hopf algebroids. Such pairings will be considered in the next section.
\ssnl
In this section, we focus on pairings of weak multiplier Hopf algebras. If we start with a regular weak multiplier Hopf algebra $A$ with enough integrals and if we take the dual $\widehat A$, we will get a special case of such a pairing.
\ssnl
We do not plan to study pairings of weak multiplier Hopf algebras in detail, but we do obtain some interesting properties. They are valuable for the discussion on pairings of algebroids we will have in the next section. A detailed and more fundamental investigation of pairings of weak multiplier Hopf algebras, such as we have for multiplier Hopf algebras in \cite{Dr-VD}, is desirable, but not necessary here. We refer to the last section, where we draw some conclusions and discuss further research for this.

\nl
\bf Dual pairs of weak multiplier Hopf algebras \rm
\nl
Let $A$ and $A'$ be two non-degenerate algebras.  We assume  that $(a,a')\to \langle a,a'\rangle$ is a non-degenerate bilinear form on $A\times A'$. Throughout, we will assume  the following basic property of such a  pairing.

\defin\label{defin:3.1}
Assume that there exist unital left and unital right actions of one algebra on the other defined by the equalities
\begin{align*}
\langle a,a'b' \rangle &= \langle a\tl a',b' \rangle  &\langle a,a'b' \rangle &=\langle  b'\tr a, a'\rangle\\
\langle ab,a' \rangle &= \langle a,b\tr a' \rangle &\langle ab,a' \rangle &=\langle  b, a'\tl a\rangle
\end{align*}
for $a,b\in A$ and $a',b'\in A'$. Then we call the pairing of the algebras \emph{admissible}.
\edefin

That the right action $\tl$ of $A'$ on $A$ is \emph{unital} means that all elements of $A$ are a linear combination of elements of the form $a\tl a'$ where $a\in A$ and $a'\in A'$. Similarly for the other actions.
\ssnl
These assumptions are quite natural in the theory.  They are  assumed e.g.\ for a pairing of multiplier Hopf algebras, see Lemma 2.3 of \cite{Dr-VD}. Moreover, it turns out that also for the examples we consider further, this condition is  satisfied.
The conditions allow us to extend the pairing as follows.
\prop
Assume that $A$ has local units. Then there is a pairing of the multiplier algebra $M(A)$ with $A'$ satisfying and characterized by
\begin{equation*}
\langle am,a'\rangle= \langle m,a'\tl a\rangle
\tussenen
\langle ma,a'\rangle=\langle m,a\tr a' \rangle
\end{equation*}
whenever $a\in A$, $m\in M(A)$ and $a'\in A'$. The pairing extends the original pairing between $A$ and $A'$.
\eprop
\bew
i) If $m$ is an element of $A$, we have that $\langle m,a'\tl a\rangle=\langle am,a'\rangle$ for all $a\in A$ and $a'\in A'$.  We will use this formula to define $\langle m,a'\rangle$ for any $m\in M(A)$ and $a'\in A'$.
\ssnl
ii) Assume that $a_i$ and $a'_i$ are finitely many elements in $A$ and $A'$ resp. If $\sum_i a'_i\tl a_i=0$ then for all $b\in A$ and $m\in M(A)$ we have
\begin{equation*}
\sum_i \langle a_imb,a'_i\rangle=\sum_i\langle mb,a'_i\tl a_i\rangle=0.
\end{equation*}
Because $A$ has local units, given $a_i$ and $m$ there is an element $b\in A$ satisfying $a_imb=a_im$ for all $i$. It follows that $\langle a_im,a'_i\rangle=0$. Because the left action of $A'$ on $A$ is assumed to be unital, all elements $a'$ in $A$ are of the form $\sum_i a'_i\tl a_i$. This shows that we can define
\begin{equation*}
\langle m,a'\tl a\rangle=\langle am,a' \rangle.
\end{equation*}
\ssnl
iii) Now we have, for all $a,b\in A$, $m\in M(A)$ and $a'\in A'$ that
\begin{align*}
\langle ma,a'\tl b \rangle
&=\langle bma,a'\rangle \\
&=\langle bm,a\tr a' \rangle\\
&=\langle m, a\tr (a'\tl b)\rangle.
\end{align*}
Because the rigth action of $A$ on $A'$ is unital we get also $\langle ma,a'\rangle=\langle m,a\tr a'\rangle$ for all $a\in A$, $m\in M(A)$ and $a'\in A'$. This completes the proof.
\ebew
We have a similar result on the other side. This means that we can define $\langle a,m'\rangle$ for $a\in A$ and $m'\in M(A')$ satisfying
\begin{equation*}
\langle a,m'a'\rangle=\langle a'\tr a,m'\rangle
\tussenen
\langle a,a'm'\rangle=\langle a\tl a',m'\rangle
\end{equation*}
whenever $a\in A$, $a'\in A'$ and $m'\in M(A')$. We can do the same for the natural tensor product pairing of $A\ot A$ with $A'\ot A'$.

\opm
A similar result is obtained in Proposition 2.3 and Proposition 2.4 of \cite{VD-W6}. There  two cases are treated in a different way. For one case we use the existence of local units while for the other case we use the existence of a counit. In this paper, we are going to apply this for weak multiplier Hopf algebras. It is shown in Propostion 2.21 of \cite{VD-W5} that local units exist for weak multiplier Hopf algebras.
\eopm
\ssnl
For finite-dimensional algebras, the product in one algebra induces a coproduct on the other one. For the general case we consider here, the existence of these coproducts is no longer true. In the first place, we can define $\Delta(a)$ for $a\in A$ as a linear functional on $A'\ot A'$ by $\langle\Delta(a),a'\ot b'\rangle=\langle a,a'b' \rangle$. Provided the algebras have local units so that the pairings can be extended, it makes sense to require that $\Delta(a)\in M(A\ot A)$. Clearly $\Delta$ is completely determined by the pairing if it exists like this.
\ssnl
This observation leads us to the following definition of a pairing of weak multiplier Hopf algebras.

\defin\label{defin:3.4}
Let $(A,\Delta)$ and $(A',\Delta')$ be weak multiplier Hopf algebras. Let $\langle \,\cdot\, ,\,\cdot\, \rangle$  be an admissible pairing of the algebras $A$ and $A'$. We call this a pairing of weak multiplier Hopf algebras if
\begin{equation}
\langle a\ot b,\Delta'(a')\rangle=\langle ab,a'\rangle
\tussenen
\langle \Delta(a),a'\ot b'\rangle=\langle a,a'b'\rangle\label{eqn:3.1}
\end{equation}
for all $a,b\in A$ and $a',b'\in A'$ and if $\langle S(a),a' \rangle =\langle a, S'(a') \rangle$ for all $a\in A$ and $a'\in A'$ where $S$ is the antipode of $A$ and $S'$ the antipode of $A'$.
\edefin

We know that the underlying algebra of a weak multiplier Hopf algebra has local units (cf. Proposition 2.21 in \cite{VD-W5}) so that the extended pairings exist and the above conditions make sense.
\ssnl
It is now easy to obtain the following consequences. We assume in what follows that we have a pairing of weak multiplier Hopf algebras as in the above definition.
\ssnl
We denote by $T_1,T_2$ and $T'_1,T'_2$ the canonical maps of $(A,\Delta)$ and $(A',\Delta')$ respectively. Recall the formulas from Equation (\ref{eqn:1.1}).

\prop\label{prop:3.5}
We have, for all $a,b\in A$ and $a',b'\in A'$, that
\begin{align}
\langle T_1(a\ot b) , a'\ot b' \rangle &= \langle a\ot b,T'_2(a'\ot b')\rangle \label{eqn:3.2}\\
\langle T_2(a\ot b) , a'\ot b' \rangle &= \langle a\ot b,T'_1(a'\ot b')\rangle \label{eqn:3.3}
\end{align}
\eprop

\bew
i) Let $a,b\in A$ and $a',b'\in A'$. Then
\begin{align*}
\langle T_1(a\ot b) , a'\ot b' \rangle
&=\langle \Delta(a)(1\ot b),a'\ot b'\rangle \\
&=\langle \Delta(a),a'\ot (b\tr b')\rangle \\
&=\langle a,a'(b\tr b')\rangle=\langle a\tl a',b\tr b'\rangle.
\end{align*}
A similar calculation will give that also $\langle a\ot b,T'_2(a'\ot b')\rangle=\langle a\tl a',b\tr b'\rangle$.
\ssnl
ii) For the two other expressions in Equation (\ref{eqn:3.3} we get $\langle b'\tr b,a'\tl a\rangle$.
\ebew

\opm
In Definition \ref{defin:3.4} we use the exented pairings to express the property that the coproducts on one algebra are adjoint to the products in the other one. The result in Proposition \ref{prop:3.5} is another way to express this property. However, as it was pointed out by the referee on an earlier version of this manuscript, it is not clear that the result of Proposition \ref{prop:3.5} will be sufficient to guarantee that also the formulas in Definition \ref{defin:3.4} hold.
\eopm

If the weak multiplier Hopf algebras are regular, we also have the two other canonical maps $T_3,T_4$ for $A$ and $T'_3,T'_4$ for $A'$, see Equation (\ref{eqn:1.2}). The corresponding equalities for these maps are
\begin{align}
\langle T_3(a\ot b) , a'\ot b' \rangle &= \langle b\ot a,T'_3(b'\ot a')\rangle\label{eqn:3.4}\\
\langle T_4(a\ot b) , a'\ot b' \rangle &= \langle b\ot a,T'_4(b'\ot a')\rangle\label{eqn:3.5}
\end{align}
for all $a,b\in A$ and $a',b'\in A'$. They again follow from the Equations (\ref{eqn:3.1}). Observe the presence of the flip maps in the right hand sides of these equalities. This is a consequence of the choices made to define these canonical maps.
\ssnl

Given the antipode, the generalized inverses of the canonical maps $T_1,T_2$ can be expressed in terms of the canonical maps $T_3,T_4$. Therefore, it should not be a surprise that we have formulas, similar to Equations (\ref{eqn:3.2}) and (\ref{eqn:3.3}), for the generalized inverses of the canonical maps.

\prop\label{prop:3.7}
Assume that we have a pairing of regular weak multiplier Hopf algebras as in Definition \ref{defin:3.4}. Then also
\begin{align}
\langle R_1(a\ot b) , a'\ot b' \rangle &= \langle a\ot b,R'_2(a'\ot b')\rangle \label{eqn:3.6}\\
\langle R_2(a\ot b) , a'\ot b' \rangle &= \langle a\ot b,R'_1(a'\ot b')\rangle \label{eqn:3.7}
\end{align}
for all $a,b\in A$ and $a',b'\in A'$. We consider here the generalized inverses of the four canonical maps $T_1,T_2$ an $T'_1,T'_2$ respectively.
\eprop

\bew
Let $a,b\in A$. We have, using the Sweedler notation,
\begin{align*}
R_1(a\ot b)
&=\sum_{(a)}a_{(1)}\ot S(a_{(2)})b\\
&=\sum_{(a)}a_{(1)}\ot S(S\inv(b)a_{(2)})\\
&=(\iota\ot S)T_3(a\ot S\inv(b)).
\end{align*}
Then, if now also $a',b'\in A'$, we get
\begin{align*}
\langle R_1(a\ot b) , a'\ot b' \rangle
&=\langle (\iota\ot S)T_3(a\ot S\inv(b)), a'\ot b' \rangle\\
&=\langle T_3(a\ot S\inv(b)), a'\ot S(b') \rangle\\
&=\langle S\inv(b)\ot a, T'_3(S(b')\ot a') \rangle.
\end{align*}
Next, using the Sweedler notation for $\Delta'$, we get
\begin{equation*}
T'_3(S(b')\ot a')=(1\ot a')\Delta'(S(b'))=\sum_{(b')}S(b'_{(2)})\ot a' S(b'_{(1)}).
\end{equation*}
If we insert this above we finally get
\begin{align*}
\langle R_1(a\ot b) , a'\ot b' \rangle
&=\sum_{(b')}\langle S\inv(b)\ot a, S(b'_{(2)})\ot a' S(b'_{(1)})\rangle\\
&=\sum_{(b')}\langle b\ot a, b'_{(2)}\ot a' S(b'_{(1)})\rangle\\
&=\sum_{(b')}\langle a\ot b, a' S(b'_{(1)})\ot b'_{(2)}\rangle \\
&=\langle a\ot b,R'_2(a'\ot b')\rangle.
\end{align*}
The other equation is proven in a similar way.
\ebew

The use of the Sweedler notation is justified as we can provide the necessary coverings. For formulas with a pairing, the covering is possible because we have the natural actions (Definition \ref{defin:3.1}) and these are unital.
\ssnl
There exist similar formulas for the generalized inverses of the other canonical maps, but we will not need these formulas.
\ssnl
\opm\label{opm:3.2}
i) In the case of a pairing of multiplier Hopf algebras, the condition $\langle S(a),a' \rangle =\langle a, S'(a') \rangle$ for all $a,a'$ can be obtained from the other axioms. The reason is that the antipode is obtained from the inverses of the canonical maps, and these are completely determined.
It is not clear if this is possible for a pairing of weak multiplier Hopf algebras. The antipode is still determined by the generalized inverses of the canonical maps, but they are obtained from a very specific choice. In order to get the condition $\langle S(a),a' \rangle =\langle a, S'(a') \rangle$ for all $a,a'$, one seems to  need an extra condition to express the compatibility of the choices made for the two weak multiplier Hopf algebras. A natural one would be to require the equalities obtained in Proposition\ref{prop:3.7} . These are obvious for multiplier Hopf algebras as the canonical maps are bijective. But it seems nicer to require the relation of the antipodes and prove the relations in Proposition \ref{prop:3.7}.
\ssnl
ii) We do not intend to study these pairings as such. In particular, we do not look for a possible minimal set of assumptions. This is an interesting project and we will say more about it in the last section where we discuss future research.
\ssnl
iii) \emph{In what follows}, we will only consider pairs of \emph{regular} weak multiplier Hopf algebras.
\eopm

Remark that, as we mentioned already, we need that the antipodes on $A$ and $A'$ are adjoints of each other in order to obtain these formulas for the generalized inverses.
\ssnl
Before we continue, consider the following expected result.

\prop Let $A$ be a regular weak multiplier Hopf algebra with enough integrals. Take for $A'$ the dual $\widehat A$. The pairing of $A$ with $A'$ is a pairing of weak multiplier Hopf algebras as in Definition \ref{defin:3.4}.
\eprop

The result is already found in Section 2 of \cite{VD-W6}. In Proposition 2.3 of that paper, it is shown that the actions induced by the duality exist and that they are unital. The formulas (\ref{eqn:3.2}) and (\ref{eqn:3.3}) in Definition \ref{prop:3.4} are used in \cite{VD-W6} to define the coproduct on the dual.
From Proposition 2.1 in \cite{VD-W6} we know that $\langle \Delta(a),a'\ot b'\rangle=\langle a, a'b'\rangle$ hollds for all $a\in A$ and $a',b'\in A'$. On the other hand, $\langle a\ot b,\Delta'(a')\rangle=\langle ab,a'\rangle$ for $a,b\in A$ and $a'\in A'$, is proven in Proposition 2.8 of \cite{VD-W6}. Observe that the first result is implicitly used in the proof of this proposition.
\ssnl
Finally, the formula relating the antipodes is used to define the antipode on the dual, see Proposition 2.11 in \cite{VD-W6}.
\nl
\bf The anti-isomorphisms of  the source and target algebras \rm
\nl
If we now combine Equation (\ref{eqn:3.2}) with (\ref{eqn:3.6}) and Equation (\ref{eqn:3.3}) with (\ref{eqn:3.7}) we arrive at
\begin{align}
\langle E(a\ot b),a'\ot b' \rangle &= \langle a\ot b ,(a'\ot 1)F'_2(1\ot b') \rangle\label{eqn:3.8} \\
\langle (a\ot b)E, a'\ot b'\rangle &= \langle a\ot b ,(a'\ot 1)F'_1 (1\ot b'\rangle\label{eqn:3.9}
\end{align}
where
\begin{equation*}
F'_1=(\iota\ot S')E'
\tussenen
F'_2=(S'\ot\iota)E'.
\end{equation*}
Here $E$ and $E'$ are the canonical idempotents of $A$ and $A'$ respectively. We have used that
$T_1R_1(a\ot b)=E(a\ot b)$ and $T_2R_2(a\ot b)=(a\ot b)E$, as well as the formulas for $R'_1T'_1$ and $R'_2T'_2$ (see Proposition \ref{prop:1.4}).
\ssnl
The following then is  a consequence of these formulas. We use $B$ and $C$ for the source and target algebras of $A$ (as before) and $B'$ and $C'$ for the source and target algebras of $A'$.

\prop\label{prop:3.4}
There exists an anti-isomorphism $\alpha:B\to B'$ and an anti-isomorphism $\beta:C\to C'$ satisfying
\begin{equation*}
\langle xa,a' \rangle= \langle a,a'S'(\alpha(x)) \rangle
\tussenen
\langle ya,a' \rangle= \langle a,\beta(y)a' \rangle
\end{equation*}
for all $a\in A$ and $a'\in A'$, where $x\in B$ and $y\in C$. We have $(\alpha\ot\beta)E=E'$.
\eprop

\bew
i) We use a Sweedler type notation for $E$ and $E'$. Take elements $b\in A$ and $b'\in A'$. First let
$x=E_{(1)}\langle E_{(2)}b,b'\rangle$. From Equation (\ref{eqn:3.8}) we find
\begin{equation*}
\langle xa,a' \rangle=\langle a, a'S'(p) \rangle
\end{equation*}
for all $a\in A$ and $a'\in A'$, where $p=E'_{(1)}\langle b,E'_{(2)}b'\rangle$. The element $p$ is uniquely defined in $B'$ and we put $\alpha(x)=p$. We claim that all elements $x\in B$ have this form and similarly for elements in $B'$. Then it will follow that $\alpha$ is an anti-isomorphism from $B$ to $B'$.
\ssnl
ii) To prove the claim, observe that for all $b\in A$
\begin{equation*}
E(1\ot b)=\sum_{(b)} b_{(2)}S\inv(b_{(1)})\ot b_{(3)}.
\end{equation*}
Then from the fullness of the coproduct on $A$ and using that the pairing is non-degenerate, it follows that all elements $x$ in $B$ are linear combinations of elements of the form $E_{(1)}\langle E_{(2)}b,b'\rangle$ with $b\in A$  and $b'\in A'$. Similarly, all elements in $B'$ are linear combinations of elements of the form $E'_{(1)}\langle b,E'_{(2)}b'\rangle.$
\ssnl
iii) Next take elements $a\in A$ and $a'\in A'$ and put $y=\langle E_{(1)}a,a'\rangle E_{(2)}$. Again from Equation (\ref{eqn:3.8}) we find
\begin{equation*}
\langle yb,b' \rangle=\langle b, qb' \rangle
\end{equation*}
for all $b\in A$ and $b'\in A'$, where $q=\langle a,a'S'(E'_{(1)})\rangle E_{(2)}'$. The element $q$ is again uniquely defined in $C'$ and all elements in $C$ are of the form $y$ above.  We put $\beta(y)=q$ and we find that $\beta$ is an anti-isomorphism from $C$ to $C'$.
\ssnl
iv)
It is an immediate consequence of the definitions of $\alpha$ and $\beta$ that $(\alpha\ot\beta)E=E'$.

\ebew

We have  similar formulas with elements of $B$ and $C$ on the other side.  They are expressed in terms of the same anti-isomorphisms.
\ssnl
In order to prove these formulas, we first need the following relation between $\alpha$ and $\beta$. They  are a consequence of the equality $(\alpha\ot\beta)E=E'$.

\prop\label{prop:3.7a}
We have $\alpha(x)=S'(\beta(S(x)))$ for all $x\in B$ and $\beta(y)=S'\alpha(S(y))$ for all $y\in C$.
\eprop

\bew
i) Take any element $x\in B$. On the one hand we have
\begin{equation*}
(\alpha\ot\beta)(E(x\ot 1))=(\alpha(x)\ot 1)(\alpha\ot\beta)(E)=(\alpha(x)\ot 1)E'.
\end{equation*}
On the other hand
\begin{align*}
(\alpha\ot\beta)(E(x\ot 1))
&=(\alpha\ot\beta)(E(1\ot S(x)))\\
&=(1\ot \beta(S(x)))(\alpha\ot\beta)(E)\\
&=(1\ot \beta(S(x)))E'.
\end{align*}
Hence, using the fullness of $E'$, we must have $\alpha(x)=S'(\beta(S(x)))$.
\ssnl
ii) Similarly take $y\in C$. We have
\begin{equation*}
(\alpha\ot\beta)((1\ot y)E)=(\alpha\ot\beta)(E)(1\ot \beta(y))=E'(1\ot \beta(y)),
\end{equation*}
as well as
\begin{align*}
(\alpha\ot\beta)((1\ot y) E)
&=(\alpha\ot\beta)((S(y)\ot 1)E)\\
&=(\alpha\ot\beta)(E)(\alpha(S(y))\ot 1)\\
&=E'(\alpha(S(y))\ot 1).
\end{align*}
This implies $S'\alpha(S(y)))=\beta(y)$.
\ebew

From these two equalities, we find that $\alpha(S^2(x))={S'}^{-2}(\alpha(x))$ for all $x\in B$. This should not come as a surprise. From the equality $E'=(\alpha\ot\beta)E$ it follows that $\mu_B\circ\alpha=\mu_{B'}$. The modular automorphism of $\mu_B$ is $S^{-2}$ while the one of $\mu_{B'}$ is ${S'}^{-2}$. Precisely, because $\alpha$ is an anti-isomorphism of $B$ with $B'$, the modular automorphism of  $\mu_B$ is converted to the \emph{inverse} of the modular automorphism of $\mu_{B'}$.
\ssnl
Now we can show the other type relations we announced.

\prop\label{prop:3.12}
We have
\begin{equation*}
\langle ax,a'\rangle = \langle a,a'\alpha(x)\rangle
\tussenen
\langle ay,a'\rangle = \langle a,S'(\beta(y))a'\rangle
\end{equation*}
for all $a\in A$ and $a'\in A'$, where $x\in B$ and $y\in C$.
\eprop

\bew
In the proof below, we will denote the antipode on $A'$ also by $S$ to simplify the expressions.
\ssnl
For the first formula, take $a\in A$, $a'\in A'$ and $x\in B$. Then
\begin{align*}
\langle ax,a'\rangle
&= \langle S(ax),S\inv(a')\rangle \\
&= \langle S(x)S(a),S\inv(a')\rangle \\
&= \langle S(a),\beta(S(x))S\inv(a')\rangle \\
&= \langle S(a),S\inv(\alpha(x))S\inv(a')\rangle \\
&= \langle S(a),S\inv(a'(\alpha(x))\rangle \\
&= \langle a,a'\alpha(x)\rangle.
\end{align*}
For the second formula, we again take $a\in A$, $a'\in A'$ but now $y\in C$. Then
\begin{align*}
\langle ay,a'\rangle
&= \langle S\inv(ay),S(a')\rangle \\
&= \langle S\inv(y)S\inv(a),S(a')\rangle \\
&= \langle S\inv(a),S(a')S(\alpha(S\inv(y)))\rangle \\
&= \langle S\inv(a),S(\alpha(S\inv(y))a')\rangle \\
&= \langle a,\alpha(S\inv y))a'\rangle \\
&= \langle a,S(\beta(y))a'\rangle.
\end{align*}
\vskip -10pt\ebew

We have used the formulas from the previous two propositions.
\ssnl
The result can also be obtained from Equation (\ref{eqn:3.7}), but by using the antipode as in the proof, we get expressions in terms of the anti-isomorphisms $\alpha$ and $\beta$.
\ssnl
Remark that the existence of these anti-isomorphisms is in agreement with the duality of weak multiplier Hopf algebras as studied in \cite{VD-W6}, see e.g.\ the item on the source and target algebras of the dual in Section 2, in particular Proposition 2.22 in \cite{VD-W6}.
\nl
Given regular weak multiplier Hopf algebras $A$ and $A'$, we can associate the regular multiplier Hopf algebroids $\mathfrak A$ and $\mathfrak A'$. Assume moreover that we have a pairing of $A$ with $A'$ as discussed in this section. The results we have obtained will allow us to define appropriate bilinear forms on balanced tensor products. To do this, we will apply the general ideas from Appendix \ref{s:appB} for various pairs. Eventually, we will arrive at a pair of regular multiplier Hopf algebroids in a sense that is discussed in the following section.

\section{\hspace{-17pt}. Pairing and duality of multiplier Hopf algebroids}  \label{s:dual2}

In this section we initiate the study of pairs of regular multiplier Hopf algebroids. It is not clear if our definition will be the final one, but at least it will be such that the two obvious cases satisfy the axioms. The first case is the one we obtain from a general pairing of weak multiplier Hopf algebras (as in the previous section). The second  is what we get from a regular multiplier Hopf algebroid with integrals and its dual. We treat both cases later in the next section and we use them to discuss this (preliminary) notion of a pairing of regular multiplier Hopf algebroids. We will then also discuss the situation that is common to both cases, that is when we start with an algebraic quantum groupoid (a regular weak multiplier Hopf algebra with enough integrals) and its dual. After all, this is what motivated us to write this paper. We will explain this.
\nl
\bf Pairing of regular multiplier Hopf algebroids \rm
\nl
We will not consider pairs of arbitrary multiplier Hopf algebroids. We will assume that the underlying source and target algebras admit  faithful linear functionals with the appropriate KMS properties.  Fortunately, this property is satisfied for the two cases we are interested in.
\snl
First we take a pair of regular multiplier Hopf algebroids $\mathfrak A$ and $\mathfrak A'$, without further conditions on the underlying source and target algebras.
\ssnl
We assume a non-degenerate and admissible pairing (see Definition \ref{defin:3.1}) of the underlying algebras $A$ of $\mathfrak A$ and $A'$ of $\mathfrak A'$.
\ssnl
As we know from the previous considerations, it is quite natural to require at least the following property of this pairing.

\voorw\label{voorw:3.9a}
If $S$ is the antipode of $\mathfrak A$ and $S'$ the antipode of $\mathfrak A'$ we assume
that $\langle S(a),a'\rangle =\langle a,S'(a')\rangle$ for all $a\in A$ and $a'\in A'$.
\evoorw

We also will impose a condition that relates the source and target algebras. It is inspired again from the case of a pairing of weak multiplier Hopf algebras (see Propositions \ref{prop:3.4} and  \ref{prop:3.5}).

\voorw\label{voorw:3.14}
As before let $B$ and $C$ be the source and target algebras of $\mathfrak A$. Let also $B'$ and $C'$ denote the source and target algebras of $\mathfrak A'$. We assume that there is an anti-isomorphism $\alpha:B\to B'$  satisfying
\begin{equation}
\langle ax,a' \rangle = \langle a,a'\alpha(x) \rangle
\tussenen
\langle xa,a' \rangle = \langle a,a'S'(\alpha(x)) \rangle\label{eqn:4.1}
\end{equation}
for all $a\in A$, $a'\in A'$, $x\in B$. We also assume that there is an anti-isomorphism  $\beta:C\to C'$ satisfying
\begin{equation}
\langle ya,a' \rangle = \langle a,\beta(y)a' \rangle
\tussenen
\langle ay,a' \rangle = \langle a,S'(\beta(y))a' \rangle\label{eqn:4.2}
\end{equation}
for all $a\in A$, $a'\in A'$ and $y\in C$.
\evoorw

It is automatic that the maps $\alpha$ and $\beta$ are anti-homomorphisms if they exist. That follows from the equations. They also must be injective because the pairing is non-degenerate. To require that they are surjective is natural from the point of view of symmetry.
\ssnl
We can now also prove the following property, relating $\alpha$ and $\beta$. We have this formula  in the previous section, see Proposition \ref{prop:3.7a}. Here we reverse the argument. In the previous section, we used the formulas of Proposition \ref{prop:3.7a} to prove the formulas in Proposition \ref{prop:3.5}. Here, we have these formulas by assumption and prove the formulas of Proposition \ref{prop:3.7a} from them.

\prop
For all $x\in B$ we have $S'(\beta(S(x)))=\alpha(x)$ while $S'(\alpha(S(y)))=\beta(y)$ for all $y\in C$.
\eprop
\bew
Given $a,a'$ and $x$ we have
\begin{align*}
\langle ax , a' \rangle
&=\langle S(x)S(a) , {S'}\inv(a')\rangle\\
&=\langle S(a) ,\beta(S(x)){S'}\inv(a') \rangle\\
&=\langle a, a'S'(\beta(S(x)))\rangle
\end{align*}
so that $S'(\beta(S(x)))=\alpha(x)$. Similarly, if we start with $$\langle ax,a'\rangle=\langle S\inv(x)S\inv(a),S'(a')\rangle,$$ we find
${S'}\inv(\beta(S\inv(x)))=\alpha(x)$ for all $x\in B$. Hence $S'(\alpha(S(y)))=\beta(y)$
\ebew

As a consequence we get that ${S'}^2(\alpha(S^2(x)))=\alpha(x)$ for all $x\in B$ as well as ${S'}^2(\beta(S^2(y)))=\beta(y)$ for all $y\in C$.

\opm
The formulas in Assumption \ref{voorw:3.14} are motivated by the results for a pairing of weak multiplier Hopf algebras as proven in Proposition \ref{prop:3.4} and Proposition \ref{prop:3.5}. There they follow from the basic assumption on the pairing of weak multiplier Hopf algebras as given in Definition \ref{defin:3.1}. We will come back to this when we argue in the next section that a pairing of weak multiplier Hopf algebras (as in Section \ref{s:dual1})  gives rise to a pairing of multiplier Hopf algebroids as developed here.
\eopm

Next we impose a condition on the source and target algebras as mentioned in the beginning. It is a condition that we seem to need in order to find appropriate pairings on the various balanced tensor products involved in the theory of multiplier Hopf algebroids. Fortunately, we will see that this property is fulfilled for the two main cases we consider later.

\voorw\label{voorw:4.5}
 We assume that there is a faithful linear functional $\mu_B$ on $B$ satisfying $\mu_B(x_1x)=\mu_B(S^2(x)x_1)$ for all $x,x_1$ in $B$.
\evoorw

The composition $\mu_B\circ S$ will be a faithful linear functional $\mu_C$ on $C$ satisfying $\mu_C(yy_1)=\mu_C(y_1S^2(y))$ for all $y,y_1$ in $C$. In what follows we will fix a choice for $\mu_B$ on $B$ and then take $\mu_C=\mu_B\circ S$. Remark that $\mu_B=\mu_C\circ S$ because in fact both functionals are invariant under the square of the antipode.
\ssnl
The formulas in Assumption \ref{voorw:3.14} guarantee that $B'$ and $C'$ are anti-isomorphic with $B$ and $C$ respectively. This implies that they have the same property. We take $\mu_{B'}$ and $\mu_{C'}$ so that
\begin{equation*}
\mu_B(x)=\mu_{B'}(\alpha(x))
\tussenen
\mu_C(y)=\mu_{C'}(\beta(y))
\end{equation*}
for all $x\in B$ and $y\in C$. Because $\alpha$ is an anti-isomorphism and because $\alpha(S^2(x)={S'}^{-2}(\alpha(x))$ we will still have that ${S'}^{-2}$ is the modular automorphism for $\mu_{B'}$. We will still have that  $\mu_{C'}=\mu_{B'}\circ S'$ and again ${S'}^2$ will be the modular automorphism of $\mu_{C'}$.
\ssnl

Remark once more that this is \emph{not a condition on the pairing}, but rather  one on each of the components. And as we see above, combined with the Assumption \ref{voorw:3.14}, it is enough to impose the condition on one component. It then follows automatically for the other one. We refer also to Remark \ref{opm:4.11a}.
\ssnl
In what follows, we make these choices for the linear functionals $\mu_B$ and $\mu_C$ on $B$ and $C$ respectively, as well for the ones on $B'$ and $C'$.
\ssnl
The next condition  formulates a restriction on the linear functionals on $A$, coming from elements in $A'$ in terms of the algebra $B$ and its faithful functional $\mu_B$.

\voorw\label{voorw:3.8}
Consider the source algebra $B$ of $\mathfrak A$ with its linear functional $\mu_B$.  We assume that for all $a \in A$  and $a'\in A'$ there exist elements $x_1$ and $x_2$ in $B$ with the property that
\begin{equation*}
\langle xa,a' \rangle = \mu_B(xx_1)
\tussenen
\langle ax,a' \rangle =\mu_B(x_2x)
\end{equation*}
for all $x\in B$.
\evoorw

We can formulate an equivalent condition in terms of the functional $\mu_C$ by using that the antipode converts $\mu_B$ to $\mu_C$.
\ssnl
The assumptions mean that the linear maps
\begin{equation*}
x\mapsto \langle xa,a' \rangle
\tussenen
x\mapsto \langle ax,a' \rangle
\end{equation*}
 belong to the subspace of the dual of $B$ consisting  of \emph{reduced} linear functionals w.r.t.\ the linear functional $\mu_B$.
\ssnl
These elements $x_1$ and $x_2$ are unique if they exist because $\mu_B$ is faithful on $B$. Therefore, if the above assumption is satisfied, we have the following property.

\prop\label{prop:4.3}
There exist bilinear maps $\gamma_1,\gamma_2$ from $A\times A'$ to $B$ satisfying
\begin{equation*}
\langle xa,a' \rangle = \mu_B(x\gamma_1(a,a'))
\tussenen
\langle ax,a' \rangle =\mu_B(\gamma_2(a,a')x)
\end{equation*}
for all $a\in A$, $a'\in A'$ and $x\in B$.
\eprop

We have the obvious module properties
\begin{equation*}
\gamma_1(xa,a')=x\gamma_1(a,a')
\tussenen
\gamma_2(ax,a')=\gamma_2(a,a')x
\end{equation*}
for all $a,a'$ in $A$ and $x$ in $B$. Using that $BA=A=AB$, it is easy to see that the formulas still hold for $x\in M(B)$. With $x=1$ we then find
\begin{equation*}
\mu_B(\gamma_1(a,a'))=\langle a,a'\rangle
\tussenen
\mu_B(\gamma_2(a,a'))=\langle a,a'\rangle
\end{equation*}
for all $a,a'$. In fact these bilinear maps are completely determined by these formulas when the module properties are assumed.
\ssnl
Remark that the condition is only a property of the linear functionals on $A$ obtained from elements in $A'$.
In the terminology of \cite{T1,T2}, it means that these functionals on $A$  are \emph{factorisable}. See Section 3.4 in \cite{T1} and Section 2.3 in \cite{T2}.
\ssnl
Playing around with the basic formulas in Proposition \ref{prop:4.3}, in particular by using the first condition on the antipodes (Assumption \ref{voorw:3.9a}), one easily gets similar bilinear maps from $A\times A'$ to $C$.

\prop\label{prop:4.4}
There exist bilinear maps $\rho_1,\rho_2$ from $A\times A'$ to $C$ satisfying
\begin{equation*}
\langle ya,a' \rangle = \mu_C(y\rho_1(a,a'))
\tussenen
\langle ay,a' \rangle =\mu_C(\rho_2(a,a')y)
\end{equation*}
for all $a\in A$, $a'\in A'$ and $y\in C$.
\eprop

\bew
We denote the antipode of $A'$ also with $S$ here.
\ssnl
Given $a\in A, a'\in A'$ and $y\in C$ we have
\begin{align*}
\langle ya,a' \rangle
&=\langle S(ya),S\inv(a') \rangle \\
&=\langle S(a)S(y),S\inv(a') \rangle \\
&=\mu_B(\gamma_2( S(a),S\inv(a'))S(y)) \\
&=\mu_C(y\rho_1(a,a'))
\end{align*}
when we put $\rho_1(a,a')=S\inv(\gamma_2( S(a),S\inv(a')))$.
\ssnl
Similarly, with $\rho_2(a,a')=S\inv(\gamma_1( S(a),S\inv(a')))$ we will have
\begin{equation*}
\langle ay,a' \rangle=\mu_C(\rho_2(a,a')y)
\end{equation*}
for all $a,a',y$.
\ebew

Again these maps are characterized by their module properties
\begin{equation*}
\rho_1(ya,a')=y\rho_1(a,a')
\tussenen
\rho_2(ay,a')=\rho_2(a,a')y
\end{equation*}
and $\mu_C(\rho_1(a,a'))=\mu_C(\rho_2(a,a'))=\langle a,a'\rangle$.
\ssnl
It is clear that conversely, the existence of the maps $\rho_1$ and $\rho_2$ will imply the existence of the maps $\gamma_1$ and $\gamma_2$ in Proposition \ref{prop:4.3}. In particular they imply Assumption \ref{voorw:3.8}.
\ssnl
We need a similar property for the functionals on $A'$ given by elements of $A$. We also want a relation of the two results with each other. The key for this is Assumption \ref{voorw:3.14}. Indeed, using the formulas from that assumption,  it is possible to find the analogues of  Propositions \ref{prop:4.3} and \ref{prop:4.4}. The argument is straightforward. However, we will not need the explicit results. Instead, we describe in the following propositions the module properties of the bilinear forms with respect to the second variable. This is in fact an equivalent result but that is how we will need it.

\prop\label{prop:4.9a}
For all $a\in A$, $a'\in A'$ and $x\in B$ we have
\begin{align*}
\gamma_1(a,a'S'(\alpha(x)))&=\gamma_1(a,a')S^{-2}(x),\\
\gamma_2(a,a'\alpha(x))&=S^2(x)\gamma_2(a,a').
\end{align*}
\eprop

\bew
i) Take $x_1\in B$ and use the second formula in Equation (\ref{eqn:4.1}) to obtain
\begin{align*}
\mu_B(x_1\gamma_1(a,a'S'(\alpha(x))))
&=\langle x_1a,a'S'(\alpha(x)) \rangle\\
&= \langle xx_1a,a'\rangle\\
&=\mu_B(xx_1\gamma(a,a'))\\
&=\mu_B(x_1\gamma(a,a')S^{-2}(x)).
\end{align*}
This proves the first equation.
\ssnl
ii) Again take $x_1\in B$ and now use the first formula in Equation (\ref{eqn:4.1}. This gives
\begin{align*}
\mu_B(\gamma_2(a,a'\alpha(x))x_1)
&=\langle ax_1, a'\alpha(x)\rangle\\
&=\langle ax_1x,a'\rangle \\
&=\mu_B(\gamma_2(a,a')x_1x)\\
&=\mu_B(S^2(x)\gamma_2(a,a')x_1).
\end{align*}
This proves the second equality.
\ebew

So we have used the two formulas of Equation (\ref{eqn:4.1}) in Assumption \ref{voorw:3.14}. Similarly, the two formulas in Equation (\ref{eqn:4.2}) will give rise to the following module properties of the maps $\rho_1$ and $\rho_2$ that we have in Proposition \ref{prop:4.4}.

\prop\label{prop:4.10a}
For all $a\in A$, $a'\in A'$ and $y\in C$ we have
\begin{align*}
\rho_1(a,\beta(y)a')&=\rho_1(a,a')S^2(y),\\
\rho_2(a,S'(\beta(y))a')&=S^{-2}(y)\rho_2(a,a').
\end{align*}
\eprop

The maps $\gamma$ and $\rho$ exists in the two cases we consider in Section \ref{s:dual3}.

\opm\label{opm:4.11a}
The existence of these maps is some kind of regularity of the pairing. Observe that these properties are automatic when the base algebras are finite-dimensional. Indeed, if e.g.\ $B$ is finite-dimensional, and if $\mu_B$ is a faithful linear functional on $B$, it automatically has the KMS property. The modular automorphism may be different from the inverse of $S^2$ but that can be achieved by modifying the linear functional. Now any linear functional $\omega$ on $B$ will be reduced in the sense that there are elements $x_1,x_2$ in $B$ so that $\omega(x)=\mu_B(xx_1)=\mu_B(x_2x)$ for all $x$ in $B$.
\eopm
\nl
We now have enough information to construct the pairings on the level of the balanced tensor products.
\nl
\bf The pairings of the balanced tensor products \rm
\nl
Under all these assumptions, we can now define the appropriate pairings of the various balanced tensor products.
\ssnl
The formulas we use further are inspired, not only by the case of a pairing of weak multiplier Hopf algebras, but also by the case of a multiplier Hopf algebroid with integrals. See Section \ref{s:dual3} where we apply this to the two cases.
\snl
We need to consider four cases in the first place. In each of these cases, we will need a bilinear map $\Gamma$ on $(A\ot A)\times (A'\ot A')$ with particular properties. Here are the definitions.

\defin\label{defin:4.11}
Using the notations from Propositions \ref{prop:4.3} and \ref{prop:4.4} we define, for $a,b\in A$ and $a',b'\in A'$,
\begin{align*}
\Gamma_1(a\ot b,a'\ot b')&=\mu_C(S(\gamma_1(a,a'))\rho_1(b,b'))\\
\Gamma_2(a\ot b,a'\ot b')&=\mu_B(\gamma_2(a,a')\gamma_1(b,b'))\\
\Gamma_3(a\ot b,a'\ot b')&=\mu_B(\gamma_2(a,a')S(\rho_2(b,b')))\\
\Gamma_4(a\ot b,a'\ot b')&=\mu_C(\rho_2(a,a')\rho_1(b,b')).
\end{align*}
\edefin

We first have the more evident module properties.

\prop\label{prop:4.12}
 Let $a,b\in A$ and $a',b'\in A'$. Then we have
\begin{align*}
\Gamma_1(a\ot yb,a'\ot b')&=\Gamma_1(S\inv(y)a\ot b,a'\ot b')\\
\Gamma_2(ax\ot b,a'\ot b')&=\Gamma_2(a\ot xb,a'\ot b')\\
\Gamma_3(ax\ot b,a'\ot b')&=\Gamma_3(a\ot bS\inv(x),a'\ot b')\\
\Gamma_4(ay\ot b,a'\ot b')&=\Gamma_4(a\ot yb,a'\ot b')
\end{align*}
when $x\in B$ and $y\in C$.
\eprop

\bew
The proof is rather straightforward. Consider e.g.\ the first formula. We have
\begin{align*}
\Gamma_1(a\ot yb,a'\ot b')
&=\mu_C(S(\gamma_1(a,a'))\rho_1(yb,b'))\\
&=\mu_C(S(\gamma_1(a,a'))y\rho_1(b,b'))\\
&=\mu_C(S(S\inv(y)\gamma_1(a,a'))\rho_1(b,b'))\\
&=\mu_C(S(\gamma_1(S\inv(y)a,a'))\rho_1(yb,b'))\\
&=\Gamma_1(S\inv(y)a\ot b,a'\ot b').
\end{align*}
The others are proven in a similar way using the module properties of the maps $\gamma$ and $\rho$ as obtained  directly from the definitions in Propositions \ref{prop:4.3} and \ref{prop:4.4} .
\ebew

Now we have the less evident module properties. We will make a comment after the proof.

\prop\label{prop:4.13}
Let $a,b\in A$ and $a',b'\in A'$. Then we have
\begin{align*}
\Gamma_1(a\ot b,a'y'\ot b')&=\Gamma_1(a\ot b,a'\ot y'b')\\
\Gamma_2(a\ot b,a'x'\ot b')&=\Gamma_2(a\ot b,a'\ot b'S\inv(x'))\\
\Gamma_3(a\ot b,a'x'\ot b')&=\Gamma_3(a\ot b,a'\ot x'b')\\
\Gamma_4(a\ot b,a'\ot y'b')&=\Gamma_4(a\ot b,S\inv(y')a'\ot b')
\end{align*}
when $x'\in B'$ and $y'\in C'$.
\eprop

\bew
Again we only prove the first formula. Take $y\in C$ and let $y'=\beta(y)$. We have
\begin{align*}
\Gamma_1(a\ot b,a'y'\ot b')
&=\mu_C(S(\gamma_1(a,a'y'))\rho_1(b,b'))\\
&=\mu_C(S(\gamma_1(a,a'\beta(y)))\rho_1(b,b'))\\
&=\mu_C(S(\gamma_1(a,a'S'(\alpha(S(y)))))\rho_1(b,b'))\\
&=\mu_C(S(\gamma_1(a,a')\sigma_B(S(y)))\rho_1(b,b'))\\
&=\mu_C(S(\gamma_1(a,a')S\inv(y))\rho_1(b,b'))\\
&=\mu_C(yS(\gamma_1(a,a'))\rho_1(b,b'))\\
&=\mu_C(S(\gamma_1(a,a'))\rho_1(b,b')\sigma_C(y))\\
&=\mu_C(S(\gamma_1(a,a'))\rho_1(b,y'b))\\
&=\Gamma_1(a\ot b,a'\ot y'b').
\end{align*}
The proof of the other cases is similar and based on the module properties of the maps $\gamma$ and $\rho$ as obtained in Propositions \ref{prop:4.9a} and \ref{prop:4.10a}.
\ebew

We can now apply the general rule, discussed in the appendix and define a pairing of the balanced spaces using the maps $\Gamma$ as follows.

\prop\label{prop:4.14}
There exist well-defined pairings
\begin{align*}
&\langle \,\cdot\, , \,\cdot\, \rangle_1  \quad\text{ on }\quad (A^C\ot {_C}A) \times (A'_C\ot {_C}A')\\
&\langle \,\cdot\, , \,\cdot\, \rangle_2  \quad\text{ on }\quad (A_B\ot {_B}A) \times (A'_B\ot {^B\hskip -3pt}A')\\
&\langle \,\cdot\, , \,\cdot\, \rangle_3  \quad\text{ on }\quad (A_B\ot {^B\hskip -3pt}A)\times (A'_B\ot {_B}A')\\
&\langle \,\cdot\, , \,\cdot\, \rangle_4  \quad \text{ on }\quad (A_C\ot {_C}A)\times ({A'}^C\ot {_C}A')
\end{align*}
induced by the maps $\Gamma_1, \Gamma_2,\Gamma_3$ and $\Gamma_4$ respectively.
\eprop

The idea is simple. Because
\begin{align*}
\Gamma_1(a\ot yb,a'\ot b')&=\Gamma_1(S\inv(y)a\ot b,a'\ot b')\\
\Gamma_1(a\ot b,a'y'\ot b')&=\Gamma_1(a\ot b,a'\ot y'b')
\end{align*}
for all $y\in C$ and $y'\in C'$ we can define
$\langle \overline{v},\overline{u}' \rangle_1 = \Gamma_1(v,u')$
where $\overline{v}$ is the image of the element $v\in A\ot A$ in $A^C\ot {_C}A$ under the quotient map and similarly for $\overline{u}'$ in $A'_C\ot {_C}A'$. Similarly for the other cases.
\nl
\bf The main definition \rm
\nl
Recall  that
\begin{equation*}
\mathfrak T_1:A_B\ot {_BA}\to A^C\ot {_CA}
\tussenen
\mathfrak T_2: A_C\ot {_CA}\to A_B\ot {^B\hskip -3pt A}
\end{equation*}
for the canonical maps $\mathfrak T_1$ and $\mathfrak T_2$ of $\mathfrak A$. We have similar properties for the canonical maps $\mathfrak T_1'$ and $\mathfrak T_2'$.
\ssnl
Therefore, the following definition makes sense.

\defin\label{defin:4.8}
Let $\mathfrak A$ and $\mathfrak A'$ be  regular multiplier Hopf algebroids with faithful linear functionals  on the source algebra and target algebras as in Assumption \ref{voorw:4.5}. We are given a non-degenerate admissible pairing  $\langle\,\cdot\,,\,\cdot\,\rangle$ of the underlying algebras $A$ and $A'$ so that $\langle S(a),a'\rangle=\langle a,S'(a')\rangle$ for all $a,a'$ where $S$ is the antipode of $\mathfrak A$ and $S'$ is the antipode of $\mathfrak A'$. For the source and target algebras $B$ and $C$ of $\mathfrak A$ and $B'$ and $C'$ of $\mathfrak A'$, we moreover assume that the pairing satisfies the Assumptions \ref{voorw:3.14} and \ref{voorw:3.8}. Finally we require that
\begin{equation}
\langle\mathfrak T_1 u,u'\rangle_1 = \langle u,\mathfrak T'_2 u'\rangle_2
\tussenen
\langle\mathfrak T_2 u,u'\rangle_3 = \langle u,\mathfrak T'_1 u'\rangle_4.\label{eqn:3.10}
\end{equation}
For the first equality we have $u\in A_B\ot {_BA}$ and $u'\in {A'}_C\ot {_CA'}$ while for the second we have $u\in A_C\ot {_CA}$ and $u'\in {A'}_B\ot {_BA'}$. The pairings are as constructed in Proposition \ref{prop:4.14}.
\ssnl
Then we call this a \emph{pairing of multiplier Hopf algebroids}.
\edefin

Now the equations with the inverses of these canonical maps follow. So we have
\begin{equation*}
\langle\mathfrak R_1 v,v'\rangle_2 = \langle v,\mathfrak R'_2 v'\rangle_1
\tussenen
\langle\mathfrak R_2 v,v'\rangle_4 = \langle v,\mathfrak R'_1 v'\rangle_3.
\end{equation*}
Here, in the first equality we have $v\in A^C\ot {_C}A$ and $v'\in A'_B\ot {^B\hskip -3pt}A'$ while in the second one we have $v\in A_B\ot {^B\hskip -3pt}A$ and $v'\in {A'}^C\ot {_C}A'$. Again we use the appropriate four pairings.
\ssnl
The formulas involving the other canonical maps will come for free by the use of the antipodes.

\opm
i) It is not clear if the above notion will be the final one. By assuming that we have an admissible pairing as in Definition 3.1, we assume the existence of the actions given in that definition. For a pairing of multiplier Hopf algebras, these actions are also assumed. They are quite natural when the pairing is obtained when taking the duality for algebraic quantum groups.
\ssnl
ii) The existence of such actions is expected, also for a pairing of multiplier Hopf algebroids as we will see when we show that it is satisfied for the pairing of a multiplier Hopf algebroid with its dual, in the case when integrals exist.
\eopm
Just as in the case of a pairing of weak multiplier Hopf algebras, we do not plan to study this concept further in this paper. We refer to Section \ref{s:concl}, where among other things we discuss possible future research, for more comments.
\ssnl
Instead, we show that the two special cases we mentioned already, satisfy our conditions. This is done in the next section.

\section{\hspace{-17pt}. Special cases}  \label{s:dual3} 

In this section we consider  two special cases of  pairs of regular multiplier Hopf algebroids.  The first case is the one we obtain from a general pairing of weak multiplier Hopf algebras (as in Section \ref{s:dual1}). The second  is what we get from a regular multiplier Hopf algebroids with integrals and its dual. We use these examples to discuss the (rather preliminary) notion of a pairing of regular multiplier Hopf algebroids we have in the previous section.
\snl
We will also discuss briefly the situation that is common to both cases, that is when we start with an algebraic quantum groupoid (a regular weak multiplier Hopf algebra with enough integrals) and its dual. After all, this is what motivated us to write this paper. We will explain this and we will say a bit more in the next section where we discuss some historical facts about the development of weak multiplier Hopf algebras, multiplier Hopf algebroids and their duality.
\nl
\bf The induced pairing from a pair of weak multiplier Hopf algebras \rm
\nl
For this case, we start with  regular weak multiplier Hopf algebras $A$ and $A'$ and a pairing of the two as in Section \ref{s:dual1}. We consider the regular multiplier Hopf algebroids $\mathfrak A$ and $\mathfrak A'$ associated to $A$ and $A'$ respectively as in Sections \ref{s:prel}. We will show that the given pairing of $A$ with $A'$ satisfies Assumptions
\ref{voorw:3.9a}n \ref{voorw:3.14}, \ref{voorw:4.5} and \ref{voorw:3.8}
of the previous section. Hence that we do have a pairing of $\mathfrak A$ with $\mathfrak A'$ as in Definition \ref{defin:4.8}.
\ssnl
By assumption, the pairing of $A$ with $A'$ is an admissible pairing (i.e.\ it admits unital actions) by definition (Definition \ref{defin:3.1}). It is also assumed that $\langle S(a),a'\rangle=\langle a,S'(a')\rangle$ (Definition \ref{defin:3.1}). The existence of the anti-isomorphisms $\alpha$ and $\beta$, needed for Assumption \ref{voorw:3.14} is obtained in Proposition \ref{prop:3.4} and \ref{prop:3.5}.
\ssnl
For the faithful linear functional $\mu_B$  take the one obtained in Proposition \ref{prop:1.8}. It satisfies the required KMS property.
\ssnl
Next we consider Assumption \ref{voorw:3.8} and we find the following expressions for the maps $\gamma_1$ and $\gamma_2$ of Proposition \ref{prop:4.3}. They are now given in terms of the canonical idempotent $E$ of $A$.

\prop\label{prop:5.1a}
We can define bilinear maps
$\gamma_1$ and $\gamma_2$ from $A\times A'$ to $C$ by
\begin{equation*}
\gamma_1(a,a')=S\inv(E_{(2)})\langle E_{(1)}a,a'\rangle
\tussenen
\gamma_2(a,a')=S(E_{(2)})\langle aE_{(1)},a' \rangle
\end{equation*}
for $a\in A$ and $a'\in A'$. We use the Sweedler expression $E_{(1)}\ot E_{(2)}$ for $E$. They satisfy
\begin{equation*}
\langle xa,a' \rangle = \mu_B(x\gamma_1(a,a'))
\tussenen
\langle ax,a' \rangle =\mu_B(\gamma_2(a,a')x)
\end{equation*}
for all $a\in A$, $a'\in A'$ and $x\in B$.
\eprop
\bew
i) First recall that $E(a\ot 1)\in A\ot C$ and $(a\ot 1)E\in A\ot C$ for all $a\in A$ (see Remark \ref{opm:1.8}). It follows that $\gamma_1$ and $\gamma_2$ are well-defined as  bilinear maps from $A\times A'$ to $B$.
\ssnl
ii) From the result given in Proposition \ref{prop:1.6} we get
\begin{equation*}
((\iota \ot S\inv)E)(x\ot 1)=(1\ot x)((\iota \ot S\inv)E)
\end{equation*}
for all $x\in B$. Then we have $\gamma_1(xa,a')=x\gamma_1(a,a')$. And if we apply $\mu_B$ we find, using the formulas of Proposition \ref{prop:1.8}, that $\mu_B(\gamma_1(a,a')) =\langle a,a'\rangle$. This proves that $\gamma_1$ satisfies $\langle xa,a' \rangle = \mu_B(x\gamma_1(a,a'))$ for all $a,a'$ and $x\in B$.
\ssnl
iii)
A similar argument works for $\gamma_2$.
\ebew
On the other hand, we have the existence of the maps $\rho_1$ and $\rho_2$  of Proposition \ref{prop:4.4}.

\prop\label{prop:5.2}
We can define bilinear maps
$\rho_1$ and $\rho_2$ from $A\times A'$ to $C$ by
\begin{equation*}
\rho_1(a,a')=S(E_{(1)})\langle E_{(2)}a,a'\rangle
\tussenen
\rho_2(a,a')=S\inv(E_{(1)})\langle aE_{(2)},a' \rangle
\end{equation*}
for $a\in A$ and $a'\in A'$. We again use the Sweedler expression $E_{(1)}\ot E_{(2)}$ for $E$. They satisfy
\begin{equation*}
\langle ya,a' \rangle = \mu_C(y\rho_1(a,a'))
\tussenen
\langle ay,a' \rangle =\mu_C(\rho_2(a,a')y)
\end{equation*}
for all $a\in A$, $a'\in A'$ and $y\in C$.
\eprop

The proof is similar to the one of the previous proposition.
\ssnl
One can verify the module properties of these maps w.r.t.\ the second variable as given in Propositions \ref{prop:4.9a} and \ref{prop:4.10a}.  We leave this as an exercise for the reader. We rather have a look at concrete realizations of the $\Gamma$ maps as in Definition \ref{defin:4.11}.

\prop
The map $\Gamma_1$  is given by
$\Gamma_1(u,u')=\langle Eu,u'\rangle$
for all $u\in A\ot A$ and $u'\in A'\ot A'$ while for $\Gamma_2$ we have
\begin{equation*}
\Gamma_2(a\ot b,a'\ot b')=\langle (a\ot 1)F_1(1\ot b),a'\ot b'\rangle
\end{equation*}
for all $a,b\in A$ and $a',b'\in A'$.
\ssnl
For the map $\Gamma_3$ we find $\Gamma_3(u,u')=\langle uE,u'\rangle$ for all $u,u'$ while
\begin{equation*}
\Gamma_4(a\ot b,a'\ot b')=\langle (a\ot 1)F_2(1\ot b) , a'\ot b'\rangle
\end{equation*}
for all $a,b\in A$ and $a',b'\in A'$.
\eprop

\bew
For all $a,b\in A$ and $a',b'\in A'$ we have
\begin{align*}
\Gamma_1(a\ot b, a'\ot b')
&=\mu_C(S(\gamma_1(a,a'))\rho_1(b,b'))\\
&=\langle S(\gamma_1(a,a'))b,b'\rangle\\
&=\langle (E_{(2)} \langle E_{(1)}a,a'\rangle) b,b'\rangle\\
&=\langle E_{(1)} a,a'\rangle \langle E_{(2)} b,b'\rangle.
\end{align*}
Similar arguments will work for the three other cases.
\ebew

We know that $\langle E(a\ot b),a'\ot b'\rangle=\langle a\ot b, (a'\ot 1)F'_2(1\ot b')$ (Equation (\ref{eqn:3.8} in Section \ref{s:dual1}. From these formulas, it is immediately clear that
\begin{align*}
&\Gamma_1(a\ot yb,a'\ot b')=\Gamma_1(S\inv(y)a\ot b,a'\ot b')\\
&\Gamma_1(a\ot b,a'\ot y'b')=\Gamma_1(a\ot b,a'y'\ot b')
\end{align*}
for all $y\in C$ and $y'\in C'$.  This illustrates what we have shown in general in Propositions \ref{prop:4.12} and \ref{prop:4.13}.
\ssnl
Similar arguments work for the other three cases.
\ssnl
The  map $\Gamma_1$ is used to define a pairing of $A^C\ot {_CA}$ with ${A'}_C\ot {_CA'}$. Similarly, $\Gamma_2$ is used to define a pairing of  $A_B\ot {_BA}$ with $A'_B\ot {^B\hskip -3pt A'}$. And it follows immediately from the equation
$\langle T_1 u,u' \rangle = \langle u,T_2u'\rangle $ for all $u\in A\ot A$ and $u'\in A'\ot A'$ that
\begin{equation}
\langle\mathfrak T_1 u,u'\rangle_1 = \langle u,\mathfrak T'_2 u'\rangle_2.
\end{equation}
In a completely similar way, from the fact that $\langle T_2 v,v' \rangle = \langle v,T_1v'\rangle$ for all $v\in A\ot A$ and $v'\in A'\ot A'$, we will get
\begin{equation}
\langle\mathfrak T_2 u,u'\rangle_3 = \langle u,\mathfrak T'_1 u'\rangle_4.
\end{equation}

Hence we can summarize with the following theorem.

\stel
Assume that we have a pairing of regular weak multiplier Hopf algebras $A$ and $A'$ in the sense of Definition \ref{defin:3.1}. Then this is also a pairing of the associated regular multiplier Hopf algebroids $\mathfrak A$ and $\mathfrak A'$ in the sense of Definition \ref{defin:4.8}.
\estel

This result is no surprise. After all, the notion of a pairing of regular multiplier Hopf algebroids in Definition \ref{defin:4.8} is defined in such a way that a pairing of regular multiplier weak Hopf algebras as in Definition \ref{defin:3.1} gives rise to a pairing of the associated multiplier Hopf algebroids, that is so that the above theorem holds.
\nl
\bf A measured multiplier Hopf algebroid and its dual \rm
\nl
Now we start from a regular multiplier Hopf algebroid $\mathfrak A$ and we assume that it has a single faithful integral. We require it to satisfy all the assumptions of Definition 2.2.1 of \cite{T2}. In other words, we assume that we have a \emph{measured regular multiplier Hopf algebroid} as studied in \cite{T2}.
\snl
As before we denote by $A$ the underlying total algebra of $\mathfrak A$. We also \emph{fix a faithful left integral} $\varphi$ on $A$.
\ssnl
Recall the following property.

\prop\label{prop:3.13}
 Let $\widehat A$ be the space of linear functionals on $A$ spanned by elements of the form $\varphi(\,\cdot\,c)$ where $c$ is in $A$. This space is also spanned by elements of the form $\varphi(c\,\cdot\,)$ with $c\in A$. We can also use a right integral and we will get the same set of functionals.
\eprop

The result is found in Equation 3.1.3 in Section 3 of \cite{T2}.

\opm\label{opm:5.7}
i) The same property holds for integrals on regular weak multiplier Hopf algebras. It is based on the fundamental property of integrals
\begin{equation*}
S((\iota\ot\varphi)(\Delta(a)(1\ot b)))=(\iota\ot\varphi)((1\ot a)\Delta(b))
\end{equation*}
and a similar formula for right integrals. See  \cite{VD-W6}.
\ssnl
ii) It is expected that the result will still hold when we only have a faithful set of integrals on multiplier Hopf algebroids (with some extra assumptions as in Definition 2.2.1 of \cite{T2}). The reason to expect this is that we also have properties like above, in item i), for integrals on multiplier Hopf algebroids, see Item 2.2 in \cite{T2}. See also a remark about this in Section \ref{s:concl} where we discuss further research.
\eopm

The set $\widehat A$ is made into an algebra and it is the underlying algebra of the dual $\widehat{\mathfrak A}$. See Theorem 3.2.3 of \cite{T2}. Having a faithful integral precisely means that we obtain a non-degenerate pairing.
\ssnl
We now let $\mathfrak A'$ be $\widehat{\mathfrak A}$  and we use $A'$ for $\widehat A$. We will show systematically that the pairing of $\mathfrak A$ with  $\mathfrak A'$ satisfies the requirements of Definition \ref{defin:4.8}.
\ssnl
As a first result, we argue that the pairing of $A$ with $A'$ is an \emph{admissible} pairing in the sense of Definition \ref{defin:3.1}.

\prop
The pairing of $A$ with $A'$ admits unital left and right actions of one algebra on the other as in Definition \ref{defin:3.1}.
\eprop

\bew
i) Let $a,b\in A$ and assume that $a'$ in $A'$ has the form $\varphi(\,\cdot\, c)$ for an element $c\in A$. Then we have $\langle ab,a' \rangle = \langle a,b\tr a' \rangle$ when $b\tr a' =\varphi(\,\cdot\, bc)$. Similarly, if $a'$ has the form $\varphi(c\,\cdot\,)$  we have $\langle ab,a' \rangle =\langle  b, a'\tl a\rangle$ when $a'\tl a =\varphi(ca\,\cdot\,)$. This proves the existence of the actions of $A'$ on $A$. They are unital because $A$ is idempotent.
\ssnl
ii) For the other side, we refer to the definition of the product in $A'$ as we find it in Theorem 3.2.3 in \cite{T2}, based on Lemma 3.2.2 of \cite{T2}. It says (among other things) that the actions of $A$ on $A'$ exist. Moreover, it is stated in Item (2) of Theorem 3.2.3 in \cite{T2} that these actions are idempotent, which means the same as unital in our terminology.
\ebew

That $\langle S(a),a'\rangle=\langle a,S(a')\rangle$ for all $a\in A$ and $a'\in A'$ is proven in Theorem 3.3.8 of \cite{T2}. So  Assumption \ref{voorw:3.9a} is satisfied.

\prop
The pairing satisfies Assumption \ref{voorw:3.14}.
\eprop

\bew
Consider the results in Lemma 3.1.1 and Lemma 3.2.5 of \cite{T2}. It is proven that there are bijections  $\eta_B:B\to C'$ and $\eta_C:C\to B'$ satisfying the equations
\begin{align}
\langle a,a'\eta_B(x) \rangle& =\langle xa,a' \rangle
\tussenen
\langle a,a'\eta_C(y) \rangle=\langle aS_B\inv(y),a'\rangle, \label{eqn:5.3a}\\
\langle a, \eta_C(y)a'\rangle&= \langle ay,a'\rangle
\tussenen
\langle a, \eta_B(x)a'\rangle=\langle S_C\inv(x)a,a' \rangle
 \label{eqn:5.4a}
\end{align}
Remark that in \cite{T2}, the algebras $B$ and $C$ are identified with the algebras $C'$ and $B'$ respectively, using the maps $\eta_B$ and $\eta_C$. We will not make this identification.
\ssnl
Define $\alpha:B\to B'$ and $\beta:C\to C'$ by
\begin{equation*}
\alpha(x)=\eta_C(S_B(x))
\tussenen
\beta(y)=\eta_B(S_C(y)).
\end{equation*}
Then we find, using the two second formulas in Equation (\ref{eqn:5.3a}) and (\ref{eqn:5.4a}),
\begin{equation*}
\langle ax,a'\rangle =\langle a,a'\alpha(x) \rangle
\tussenen
\langle ya,a' \rangle=\langle a, \beta(y)a'\rangle
\end{equation*}
for all $a,a'$. As we have seen before, it follows form these two equations that
\begin{equation*}
\beta(y)=S_B'(\alpha(S_C(y)))
\tussenen
\alpha(x)=S_C'(\beta(S_B(x)))
\end{equation*}
 and it will follow that
\begin{equation*}
\eta_B(x)=S_B'(\alpha(x))
\tussenen
\eta_C(y)=S_C'(\beta(y))
\end{equation*}
for all $x,y$. Then the two first formulas in Equation (\ref{eqn:5.3a}) and (\ref{eqn:5.4a}) will give
\begin{equation*}
\langle xa,a'\rangle=\langle a,aS'(\alpha(x))\rangle
\tussenen
\langle ay,a' \rangle=\langle a, S'(\beta(y))a'\rangle
\end{equation*}
for all $a,a'$ and $x,y$. This completes the proof.
\ebew

 By the definition of a measured regular multiplier Hopf algebroid, see Definition 2.2.1 in \cite{T2} there exists a faithful linear functional $\mu_B$ on $B$ satisfying $\mu_B(x_1x)=\mu_B(S^2(x)x_1)$ for all $x,x_1\in B$, see Theorem 2.2.3 in \cite{T2}. Hence also Assumption \ref{voorw:4.5} is fulfilled.

\prop\label{prop:4.7}
The pairing of $A$ with $A'$ satisfies Assumption \ref{voorw:3.8}.
\eprop

\bew
i) Take $a\in A$ and $a'\in A'$. Assume that $a'$ has the form $\varphi(\,\cdot\,c)$ where $c\in A$. By assumption (Item (4) in Definition 2.2.1 of \cite{T2}) there is a linear map $_B\Phi:A\to B$ satisfying $_B\Phi(xa)=x\, _B\Phi(a)$ and  $\mu_B(_B\Phi(a))=\varphi(a)$ for all $a\in A$ and $x\in B$.
Then
\begin{equation*}
\langle xa,a'\rangle=\varphi(xac)=\mu_B(_B\Phi(xac))=\mu_B(xx_1)
\end{equation*}
where $x_1={_B\Phi(ac)}$.
\ssnl
ii) Again take $a\in A$ and $a'\in A'$. Now assume that $a'$ had the form $\varphi(c\,\cdot\,)$ where $c\in A$. Again from Item (4)  of Definition 2.2.1 in \cite{T2}, we know that there is a map $\Phi_B:A\to B$ satisfying $\Phi_B(ax)=\Phi_B(a)x$ and $\mu_B(\Phi_B(a))=\varphi(a)$ for all $a\in A$ and $x\in B$. Then
\begin{equation*}
\langle ax,a'\rangle=\varphi(cax)=\mu_B(\Phi_B(cax))=\mu_B(x_2x)
\end{equation*}
where now $x_2=\Phi_B(ca)$.
\ebew

We are now only one step away from the proof of our main result.
\ssnl
First we observe that the pairings we have constructed for a general pair of regular multiplier Hopf algebroids in Proposition \ref{prop:4.14} are precisely the same as the pairings constructed on the balanced tensor products in Lemma 3.3.1 of \cite{T2}.
\ssnl
And hence we arrive at the following.

\stel\label{stel:5.11}
Assume that $\mathfrak A$ is a measured regular multiplier Hopf algebra  and that $\mathfrak A'$ is its dual as constructed in Theorem 3.3.6 of \cite{T2}. Then the pairing of the underlying algebras $A$ and $A'$ yields a pairing of regular multiplier Hopf algebroids in the sense of  Definition \ref{defin:4.8}.
\estel

We essentially just need to argue that the canonical maps satisfy
\begin{equation}
\langle\mathfrak T_1 u,u'\rangle_1 = \langle u,\mathfrak T'_2 u'\rangle_2
\tussenen
\langle\mathfrak T_2 u,u'\rangle_3 = \langle u,\mathfrak T'_1 u'\rangle_4.
\end{equation}

This however is shown in Proposition 3.3.3 of \cite{T2}.
\nl
\bf Duality for weak multiplier Hopf algebras with a single faithful integral \rm
\nl
Finally, we start from a regular weak multiplier Hopf algebra $A$ with a single faithful integral. Take for $A'$ the dual $\widehat A$ and consider the associated pairing of $A$ with $A'$.
\ssnl
On the one hand, this pairing gives rise to a pairing of the associated regular multiplier Hopf algebroids $\mathfrak A$ and $\mathfrak A'$ as in the first item of this section. On the other hand, we can first  consider the regular multiplier Hopf algebroid $\mathfrak A$ associated with $A$. As we have shown in Section \ref{s:from}, it will have a single faithful integral and hence, it will be a regular measured multiplier Hopf algebroid in the sense of \cite{T2}. It allows a dual $\widehat{\mathfrak A}$. And again we obtain a pairing of  multiplier Hopf algebroids according to the second item in this section.
\ssnl
It is fairly straightforward to see that these two procedures eventually result in the same pair of multiplier Hopf algebroids. Let us consider the various steps and illustrate this.
\ssnl
In the first place we have the pairing of the underlying algebras $A$ and $A'$ that is the basis for the two procedures. The properties of this pairing are the same in the two cases. We have the requirement that the pairing of the two algebras is admissible. We have that the antipodes are adjoints of each other. We have the anti-isomorphisms $\alpha:B\to B'$ and $\beta:C\to C'$. Finally we have the existence of the maps $\gamma_1$, $\gamma_2$, $\rho_1$ and $\rho_2$. All these properties are properties of the pairing of the underlying total algebras $A$ and $A'$ and their base algebras $B,C$ and $B',C'$.
\ssnl
However, in the first case, the maps $\gamma_1,\gamma_2$ and $\rho_1,\rho_2$ are given in terms of the canonical idempotent $E$ while in the second case, they are obtained in therms of the left integrals as in Proposition \ref{prop:4.7}. We verify that these are the same.

\prop
Under the assumptions formulated above, let $\gamma_1$ and $\gamma_2$ be the maps from $A\times A'$ to $B$ as in Proposition \ref{prop:5.1a}. On the other hand, let $\gamma'_1$ and $\gamma'_2$ be the maps from $A\times A'$ to $B$ as obtained in the proof of Proposition \ref{prop:4.7}. Then $\gamma_1=\gamma'_1$ and $\gamma_2=\gamma'_2$. A similar result holds for the maps $\rho_1$ and $\rho_2$ as obtained in Proposition \ref{prop:5.2} and the maps $\rho'_1$ and $\rho'_2$ from Proposition \ref{prop:4.4}.
\eprop

\bew
i) The  maps $\gamma_1,\gamma'_1$ and $\gamma_2,\gamma'_2$  satisfy
\begin{align*}
\langle xa,a' \rangle&=\mu_B(x\gamma_1(a,a'))=\mu_B(x\gamma'_1(a,a')\\
\langle xa,a' \rangle&=\mu_B(\gamma_2(a,a')x)=\mu_B(\gamma'_1(a,a')x)
\end{align*}  as we see from Proposition \ref{prop:5.1a} and  Proposition \ref{prop:4.7} resp. Then $\gamma_1=\gamma'_1$ and $\gamma_2=\gamma'_2$ because $\mu_B$ is faithful. We are using that the map $\mu_B$ is the same for the two approaches.
\ssnl
ii) A similar argument applies for the $\rho$-maps.
\ebew

As a result of all these considerations, we get the following essential result.

\stel\label{stel:5.12}
Let $A$ be a regular weak multiplier Hopf algebra with single faithful integral and let $\widehat A$ be the dual in the sense of duality for weak multiplier Hopf algebras (as in \cite{VD-W6}). The multiplier Hopf algebroid associated to $\widehat A$ as in \cite{T-VD2} is the dual, in the sense of duality for multiplier Hopf algebroids (as in \cite{T2}) of the multiplier Hopf algebroid $\mathfrak A$ associated with $A$.
\estel

Remark that it is expected that the above result will also hold when we only start with the existence of a faithful set of integrals. However, for this one first needs the duality theory of multiplier Hopf algebroids with a faithful set of integrals, see \ref{opm:5.7}.

\section{\hspace{-17pt}. Some historical consideration} \label{s:hist} 

We will describe in this section a short history of the development of weak multiplier Hopf algebras, multiplier Hopf algebroids and their duality. Contributions are mainly by the authors of this note, as well as by B.-J. Kahng.
\nl
\bf Weak multiplier Hopf algebras - The basic theory\rm
\nl
The initiative for working on weak multiplier Hopf algebras came originally from one of us (SHW). The first joint (but unpublished) work (AVD \& SHW)  was entitled \emph{Multiplier Unifying Hopf algebras} (2009). This was  the basis for our work on weak multiplier Hopf algebras.
\ssnl
The first two papers (by AVD \& SHW) written on the subject were available on the Arxiv in 2012. In  \emph{Weak multiplier Hopf algebras. Preliminaries, motivation and basic examples}  (arXiv:1210.3954) we explained our ideas and our motivation for the definition of a weak multiplier Hopf algebra. In \emph{Multiplier Hopf algebras I. The main theory} (arXiv:1210.4395) we developed the main theory. The first paper appeared in an issue of the Banach Center Publications already in 2012 \cite{VD-W3}. The second one took more time and appeared in Crelles Journal in 2015 \cite{VD-W4}.  So, the notion of a \emph{weak multiplier Hopf algebra} was developed in the period 2009 - 2012 by two of us (AVD \& SHW).
\ssnl
The \emph{source and target algebras}, together with the source and target maps, were introduced already in Section 3 of \cite{VD-W3}. A more complete study appeared on the Arxiv  in 2014 under the title \emph{Weak multiplier Hopf algebras II. The source and target algebras} (arXiv:1403.7906). A second version of this work dates from 2015 and is substantially different from the first version in the sense that it also treats the non-regular case.  This work is not (yet) published. The reference is \cite{VD-W5}.
\ssnl
Intimately related with the study of the source and target algebras is the work on \emph{separability idempotents} in the multiplier algebra setting. Also here we have a first version on the Arxiv as early as in 2013 (arXiv:1301.4398v1) where only the regular case is treated while there is a second version, again substantially different from the first one, put on the Arxiv in 2015 (arXiv:1301.4398v2), where also the non-regular case is treated. Also this work is not (yet) published. The references here  \cite{VD4.v1, VD4.v2}
\ssnl
The main theory of weak multiplier Hopf algebras  is from the very beginning developed also for the non-regular case. Roughly speaking, the regular case is such that the antipode is a bijection of the underlying algebra. There are examples of Hopf algebras with an antipode that does not have this property. It means that there are weak multiplier Hopf algebras that are not regular. However, more research is needed here in order to provide more interesting (from the point of view of weak multiplier Hopf algebras) non-regular cases. This seems to be rather difficult, but on the other hand very few attempts have been made to search such examples. See some more comments and suggestions about this in the next section where we discuss possible further research.
\ssnl
The development on this level dates from about 2013 till 2015.
\nl
\bf Integral theory and duality \rm
\nl
Almost from the very beginning, when we started the development of weak multiplier Hopf algebras, we (AVD \& SHW) started with looking at integral theory and duality. The first versions we have about this are from 2010 and contained already most of the results on integrals. Most results on the dual we had in 2011. These results were presented at several talks. The first occasion was at a conference in Warsaw (2011), see \cite{VD3}. The second was at a conference in Caen (2012). At the moment of writing this note, the slides of this talk are still available on the net, see \cite{VD3a}. Let us mention also the talk at the University of Budapest (Hungary) in 2013, see \cite{VD3b}. However, the first version of these results was only available on the Arxiv in 2017 (\cite{VD-W6}). The reason why it took so many years for this is simply because we were in the mean time working on various other projects while on the other hand, we used a fair amount of time to provide a clear and complete presentation of our results.
\ssnl
In 2014 a preprint by Byung-Jay Kahng and AVD was put on the Arxiv including results on integrals on weak multiplier Hopf algebras (arXiv:1406.0299). It was entitled \emph{The Larson Sweedler theorem for weak multiplier Hopf algebras}. This work was inspired by the work of B\"ohm et al. \cite{B-G-L} on weak multiplier bialgebras that was available already in 2014. The Larson Sweedler paper appeared in 2018 (see \cite{K-VD}).
\ssnl
As we see, the work on integrals and duality of weak multiplier Hopf algebras has a longer history, but was essentially developed in the period 2010-2014.
\nl
\bf Multiplier Hopf algebroids \rm
\nl
Soon after the development of the main theory of weak multiplier Hopf algebras, we started to work on an Hopf algebroid version. This was done by two of us (TT and AVD). The first version of the work appeared on the Arxiv in 2013 (arXiv:1307.0769). It took a while before it was actually published in \cite{T-VD1} (2017). The relation of these multiplier Hopf algebroids with weak multiplier Hopf algebras was obtained in 2014 (see  arXiv:1406.3509). The work appeared in the Banach Center Publications in 2015 \cite{T-VD2}. Related is the paper by AVD on \emph{Modified weak multiplier Hopf algebras}, unpublished but available on the Arxiv since 2014 (arXiv:1407.0513).
\ssnl
Multiplier Hopf algebroids were developed in the period 2012-2014.
\nl
\bf Integrals on multiplier Hopf algebroids and duality\rm
\nl
Next comes the theory of integrals and duality for multiplier Hopf algebroids, developed by one of us (TT) during the period 2014-2016. First there is the work on integrals on the Arxiv (arXiv:1403.5282) entitled \emph{Regular multiplier Hopf algebroids II. Integration on and duality of algebraic quantum groupoids}. An improved version of the first part on integrals appeared on the Arxiv in 2015 (arXiv:1507.00660) while results on the duality are found on the Arxiv in 2016 (arXiv:1605.06384). The last two papers have appeared in 2016 \cite{T1} and 2017 \cite{T2} respectively.
\ssnl
The work on integrals and duality of multiplier Hopf algebroids by TT was inspired by the work on integrals and duality of weak multiplier Hopf algebras, developed earlier by AVD, SHW and Kahng, but further developed independently.
\ssnl
In \cite{T2} the duality of weak multiplier Hopf algebras with integrals is obtained as an application of the duality theory of multiplier Hopf algebroids with integrals. This is without any doubt a nice illustration of the theory, but it should not be considered as a means of obtaining the duality for weak multiplier Hopf algebras as in \cite{VD-W6}. It is a too heavy machinery for this task. The direct method is a lot easier and far more natural. The theory of integrals on multiplier Hopf algebroids is nice and very interesting, but rather involved. The study of integrals and duality in \cite{VD-W6} is not only more easy, but also more concrete and also provides more information. In fact, by treating the case of a faithful set of integrals, makes it also more general.
\ssnl
With our present work, we provide the interested reader more insight in the relation of the two theories. And it should be helpful for understanding how the the algebroid theory relates with, and is inspired by the easier weak multiplier Hopf algebra duality theory.

\section{\hspace{-17pt}. Conclusion and further research}\label{s:concl}  

In this paper we initiated the study of pairings of weak multiplier Hopf algebras, of multiplier Hopf algebroids and of the relation between the two. We do not claim to have the best possible notions defined yet. This needs more research. We have an approach, mainly motivated by the basic examples. For a pairing of weak multiplier Hopf algebras, this is the one coming from a regular weak multiplier Hopf algebra with enough integrals and its dual (as constructed in \cite{VD-W6}). For a pairing of multiplier Hopf algebroids, we think of a multiplier Hopf algebroid with a single faithful integral and its dual (as in \cite{T2}). For the relation of the two notions, we are considering the passage from a pair of weak multiplier Hopf algebras to the associated pair of multiplier Hopf algebroids.
\ssnl
We believe that the notion of a pairing of weak multiplier Hopf algebras, as given in Section \ref{s:dual1} in this paper, is probably very close to the optimal notion. It is indeed inspired by the well-established version of a pairing of multiplier Hopf algebras as given in \cite{Dr-VD}. We have seen how a regular multiplier Hopf algebra with enough integral $(A,\Delta)$ is paired with its dual $(\widehat A,\widehat\Delta)$. But it would be interesting to have other examples of such a pairing. There are well-known pairings of Hopf algebras that do not come from the construction of the dual. They will fit in this more general notion, but of course, we would like to have examples of pairings of genuine multiplier Hopf algebra that are not coming from duals constructed from integrals. More research is welcome here.
\ssnl
The situation is a bit different when we look at pairings of multiplier Hopf algebroids. The notion as studied in Section \ref{s:dual2} of this paper is probably still immature. We are in this paper mainly interested in a notion sufficient to include the two cases we consider (see Section \ref{s:dual3}). But more research on this topic, both from the theoretical point of view, as what the search for examples is concerned has to be done.
\ssnl
There is the issue of regularity. Apart from the more basic results, mostly only regular weak multiplier Hopf algebras and regular multiplier Hopf algebroids are considered. Nevertheless, it seems interesting also to study the non-regular case further. There is mainly a lack of examples, apart from the ones obtained of Hopf algebras with a non-invertible antipode. On the other hand, it is not clear if an integral theory makes sense in the non-regular case. Also this should be investigated.
\ssnl
And then finally, there is the link of the algebraic theories with the operator algebraic ones. Some steps have been achieved. There is well-established notion within the von Neumann algebraic framework. However no such satisfactory theory exists in the C$^*$-algebraic setting. For locally compact quantum groups, these two frameworks provide essentially the same objects and the relation between the two is well understood (see e.g.\ \cite{VD9}). But this seems a lot more complicated in the case of quantum groupoids. There are several attempts in this direction (see e.g.\ \cite{K-VD1}) but no final and satisfactory results have been obtained yet.

\renewcommand{\thesection}{\Alph{section}}

\setcounter{section}{0}

\renewenvironment{stelling}{\begin{itemize}\item[ ]\hspace{-28pt}\bf Theorem \rm }{\end{itemize}}
\renewenvironment{propositie}{\begin{itemize}\item[ ]\hspace{-28pt}\bf Proposition \rm }{\end{itemize}}
\renewenvironment{lemma}{\begin{itemize}\item[ ]\hspace{-28pt}\bf Lemma \rm }{\end{itemize}}

\section{\hspace{-17pt}. Appendix. Adjoint linear maps}\label{s:appB} 

In this appendix, the starting point is a pair of vector spaces $X$ and $X'$ (over the field $\mathbb C$ of complex numbers) together with a non-degenerate pairing $(x,x')\mapsto \langle x,x'\rangle$ from $X\times X'$ to $\mathbb C$. Furthermore, we have a linear map $T$  from $X$ to itself and a linear map $T'$ from  $X'$ to itself satisfying $\langle Tx,x' \rangle = \langle x,T'x' \rangle$ for all $x$ in $X$ and $x'$  in $X'$. We say that $T$ and $T'$ are \emph{adjoints} of each other.
Finally, we assume that  that $T$ and $T'$ have generalized inverses $R$ and $R'$, again adjoint to each other in the sense that also $\langle Rx,x' \rangle = \langle x,R'x \rangle$ for all $x,x'$.
\ssnl
In the paper, we have several situations like this and what we do in this appendix applies to all of them (see the item on the pairing on the balanced tensor product in Sections \ref{s:dual2} and \ref{s:dual3}). 
\ssnl
The results in this appendix certainly are not very deep, in fact almost trivial, but nevertheless, they give us some extra understanding of certain aspects we encounter in this paper.

\notat
 We consider the following compositions of linear maps:
\begin{align*}
E&=TR \tussenen E'=T'R'\\
F&=RT \tussenen F'=R'T'.
\end{align*}
\enotat
Then $E$ and $E'$ are projections on the range of $T$ and $T'$ respectively, while $F$ and $F'$ are projections on the kernels. This follows from the definition of a generalized inverse. We have that $E$ is adjoint to $F'$ and $F$ adjoint to $E'$.
\ssnl
Next we consider quotient spaces and the corresponding projection maps.

\notat\label{notat:B.2}
Denote by $X_s$ and $X_t$ the kernel of $T$ and of $R$ respectively. Similarly we use $X'_s$ and $X'_t$ for the kernel of $T'$ and of $R'$. We define the quotient spaces, together with their quotient maps accordingly:
\begin{align*}
&\pi_s:X\to \overline X_s  &\pi'_s:X'\to \overline X'_s \\
&\pi_t:X\to \overline X_t  &\pi'_t:X'\to \overline X'_t
\end{align*}
where $\overline X_s=X/X_s$, $\overline X_t=X/X_t$, etc.
\enotat

We now consider the following result.

\prop\label{prop:B.3}
There are canonically associated two non-degenerate pairings
 \begin{align*}
&(v,u') \mapsto \langle v,u' \rangle_1 \quad\text{ on } \quad \overline X_t\times \overline X'_s,\\
&(u,v') \mapsto \langle u,v' \rangle_2 \quad\text{ on }\quad \overline X_s\times \overline X'_t.
\end{align*}
They are defined by
\begin{align*}
\langle \pi_t(x),\pi'_s(x') \rangle_1&:=\langle Ex,x' \rangle =\langle x,F'x' \rangle\\
\langle \pi_s(x),\pi'_t(x') \rangle_2&:=\langle Fx,x' \rangle =\langle x,E'x' \rangle
\end{align*}
\eprop
\bew
The proof is straightforward. For the first case we observe that $\pi_t(x)=0$ if and only if $Ex=0$ and that $\pi'_s(x')=0$ if and only if $F'x'=0$. Similarly for the other case.
\ebew

Next we consider the induced linear maps between the appropriate quotient spaces.

\notat
Define
\begin{align*}
&\mathfrak T: \overline X_s \to \overline X_t
\tussenen
\mathfrak T': \overline X'_s \to \overline X'_t\\
&\mathfrak R: \overline X_t \to \overline X_s
\tussenen
\mathfrak R': \overline X'_t \to \overline X'_s
\end{align*}
by $\mathfrak T\pi_s(x)=\pi_t(Tx)$ and $\mathfrak R\pi_t(x)=\pi_s(Rx)$.
\enotat
All these  maps are obviously well-defined. They are bijections, the maps $\mathfrak T$ and $\mathfrak R$ are inverses of each other and the same for the maps $\mathfrak T'$ and $\mathfrak R'$. Moreover, $\mathfrak T$ and $\mathfrak T'$ are adjoints to each other and the same for the maps $\mathfrak R$ and $\mathfrak R'$. More precisely we have:

\prop\label{prop:B.5}
We have
\begin{equation*}
\langle \mathfrak T u,u' \rangle_1= \langle u, \mathfrak T' u'\rangle_2
\tussenen
\langle \mathfrak R v,v' \rangle_2= \langle v, \mathfrak R' v'\rangle_1.
\end{equation*}
Here $u\in \overline X_s$, $u'\in \overline X'_s$ and $v\in \overline X_t$, $v'\in\overline X'_t$
\eprop
\bew
i) Consider the first formula with $u\in \overline X_s$ and $u'\in \overline X_s$. By assumption $\mathfrak Tu\in \overline X_t$ and so $\langle \mathfrak T u,u' \rangle_1$ is defined. Similarly $\mathfrak T'u'\in \overline X'_t$ and therefore $\langle \mathfrak u, Tu' \rangle_2$ is also defined.
\ssnl
ii) Now take elements $x\in X$ and $x'\in X'$ so that $\pi_s(x)=u$ and $\pi'_s(x')=u'$. Then we have
\begin{align*}
\langle \mathfrak T u,u' \rangle_1
&=\langle \mathfrak T\pi_s(x),\pi'_s(x')\rangle_1\\
&=\langle \pi_t(Tx),\pi'_s(x')\rangle_1\\
&=\langle ETx,x'\rangle=\langle Tx,x' \rangle,
\end{align*}
as well as
\begin{align*}
\langle u,\mathfrak T'u' \rangle_2
&=\langle \pi_s(x),\mathfrak T'\pi'_t(x')\rangle_2\\
&=\langle \pi_s(x),\pi'_t (T'x')\rangle_2\\
&=\langle x,E'T'x'\rangle=\langle x,T'x'\rangle.
\end{align*}
This proves the first equality.
\ssnl
iii) The proof of the other equation follows in the same manner from $\langle Rx,x'\rangle=\langle x,R'x'\rangle$.
\ebew

\opm\label{opm:B.6}
i) Suppose that we are in the situation of a pairing between vector space $X,X'$ like in this appendix. Now we consider subspaces $X_s, X_t$ of $X$ and $X'_s, X'_t$ of $X'$, but we no longer assume that they are obtained as the kernel and the range of linear maps $T$ and $T'$. The quotient spaces as in Notation \ref{notat:B.2} can be defined, together with the canonical projections.
\ssnl
ii) In general, we can not hope to define canonically the induced pairings as in Proposition \ref{prop:B.3}. However, we do not really need the linear projection maps $E,F$ and $E',F'$. It is enough to have two bilinear forms $\Gamma_1$ and $\Gamma_2$ on $X\times X$ with the property that
$\Gamma_1(x,x')=0$ when $x\in X_s$ and $x'\in X'_t$ and that $\Gamma_2(x,x')=0$ when $x\in X_t$ and $x'\in X'_s$.
\ssnl
iii) The bilinear form $\Gamma_1$ induces  linear maps $F$ from $X$ to the linear dual of $X'$ and $E'$ from $X'$ to the linear dual of $X$ so that $Fx=0$ for $x\in X_s$ and $Ex'=0$ when $x'\in X'_t$ and also $\langle Fx,x'\rangle=\langle x,E'x'\rangle$ for all $x,x'$. We have the obvious extended pairings here. Similarly, the form $\Gamma_2$ is given by linear maps $E$ and $F'$.
\ssnl
iv) The pairings of the quotient spaces are now defined with the same formulas as in Proposition \ref{prop:B.3}.
\eopm

\end{document}